\documentclass[review]{elsarticle}

\usepackage{hyperref}
\usepackage{epsfig}
\usepackage{epstopdf}
\usepackage{bbm}
\usepackage{amsmath}
\usepackage{amssymb}
\usepackage{mathrsfs}

\newcommand{\ie}{{\it i.e.}}
\newcommand{\eg}{{\it e.g.}}

\graphicspath{{./pics/}}

\journal{J. Optim. Theory Appl.}

\usepackage{numcompress}

\begin{document}

\begin{frontmatter}

\title{An Improved Multi-Parametric Programming Algorithm for Flux Balance Analysis of Metabolic Networks}

%\author{Amir Akbari and Paul I. Barton}
%\address{Process Systems Engineering Laboratory, Massachusetts Institute of Technology, Cambridge, MA 02139, United States}

\author[mymainaddress]{Amir Akbari}
\ead{aakbari@mit.edu}
\author[mymainaddress]{Paul I. Barton\corref{mycorrespondingauthor}}
\cortext[mycorrespondingauthor]{Corresponding author}
\ead{pib@mit.edu}
\address[mymainaddress]{Process Systems Engineering Laboratory, Massachusetts Institute of Technology, Cambridge, MA 02139, United States}

\begin{abstract}
Flux balance analysis has proven an effective tool for analyzing metabolic networks. In flux balance analysis, reaction rates and optimal pathways are ascertained by solving a linear program, in which the growth rate is maximized subject to mass-balance constraints. A variety of cell functions in response to environmental stimuli can be quantified using flux balance analysis by parameterizing the linear program with respect to extracellular conditions. However, for most large, genome-scale metabolic networks of practical interest, the resulting parametric problem has multiple and highly degenerate optimal solutions, which are computationally challenging to handle. An improved multi-parametric programming algorithm based on active-set methods is introduced in this paper to overcome these computational difficulties. Degeneracy and multiplicity are handled, respectively, by introducing generalized inverses and auxiliary objective functions into the formulation of the optimality conditions. These improvements are especially effective for metabolic networks because their stoichiometry matrices are generally sparse; thus, fast and efficient algorithms from sparse linear algebra can be leveraged to compute generalized inverses and null-space bases. We illustrate the application of our algorithm to flux balance analysis of metabolic networks by studying a reduced metabolic model of \emph{Corynebacterium glutamicum} and a genome-scale model of \emph{Escherichia coli}. We then demonstrate how the critical regions resulting from these studies can be associated with optimal metabolic modes and discuss the physical relevance of optimal pathways arising from various auxiliary objective functions. Achieving more than five-fold improvement in computational speed over existing multi-parametric programming tools, the proposed algorithm proves promising in handling genome-scale metabolic models.
\end{abstract}


\end{frontmatter}

%\linenumbers

\section{Introduction} \label{sec:introduction}
Multi-parametric programming (MPP) encompasses a broad class of optimization problems with nonlinear objectives and constraints where variations of optimal solutions with respect to a set of parameters are of interest. The domain of parameters is not necessarily small, so one generally seeks a global parametric solutions in a subset of the parameter space. This contrasts with sensitivity analysis where local variations in optimal solutions with respect to infinitesimal perturbations in parameters are desired \cite{Bertsimas1997, Borrelli2003, Bonnans2013}.

This paper concerns MPP for two special classes of optimization problems, namely, linear programs (LPs) and convex quadratic programs (QPs). The objective function is linear in the decision variables in LPs and quadratic in QPs, while equality and inequality constraints are linear in both. The goal in MPP is to construct explicit parametric optimal solutions by partitioning the parameter space into critical regions (CRs) (see Definition 2.3 and Eq.~(3) of \citet{Borrelli2003}), where the solutions are determined by the respective active set in each CR. We restrict our analysis to problems where optimal solutions are affinely parametrized with respect to the right-hand side of the constraints, where MPP furnishes a polyhedral partition of the parameter space. This formulation commonly arises in a wide range of practical applications, including production planning, scheduling problems, communication-network planning, model predictive control, and metabolic-network simulations \cite{Bertsimas1997, pistikopoulos2007multi, banga2008optimization}.

Several MPP algorithms for LPs have been developed in the literature using the Simplex-tableau-basis \cite{Gal1972} and active-set methods \cite{Borrelli2003}. The latter have gained more interest due to simplicity of implementation and convenient handling of degeneracies. Therefore, subsequent algorithms for QPs were mostly based on the active-set method \cite{Dua2000, Bemporad2002, Tondel2003a, Tondel2003}. However, the multiplicity of optimal primal and dual solutions\footnote{In this paper, multiplicity (degeneracy) of optimal primal and dual solutions are, respectively, referred to as `primal multiplicity (degeneracy)' and `dual multiplicity (degeneracy)'.} can significantly complicate the construction of CRs and optimal solutions, respectively.

\citet{Mangasarian1987} showed that the continuation of optimal solutions is Lipschitz continuous with respect to the right-hand-side parameters for LPs. This implies that, irrespective of whether the optimal solution is unique or non-unique, a continuous parametric solution can always be constructed along any path in the parameter space, especially those crossing a boundary between two neighboring CRs\footnote{Note that crossing a boundary between two CRs corresponds to a change in the active set.}. When multiple optimal solutions exist, although all the solutions are acceptable with respect to the objective function, any arbitrary choice of solution, as suggested by the algorithm of \citet{Borrelli2003}, does not guarantee a continuous continuation of solutions between neighboring CRs. Moreover, different choices can lead to significantly different partitions. Therefore, it is important to establish an appropriate strategy to associate a unique solution with each point of the parameter space. Reformulating the problem as a hierarchical optimization using auxiliary objective functions is an expedient technique to identify the ``best'' among all possible solutions. Here, the properties of a specific problem determines the criteria for what the best solution is. For example, the continuity of all decision variables is highly desirable in model predictive control \cite{Bemporad2002a}, which can be ensured by strictly convex auxiliary quadratic objective functions \cite{Nocedal2006}.

To address the forgoing non-uniqueness issue, \citet{Spjotvold2007} proposed an auxiliary norm-minimizing objective function for LPs. Since this objective function is strictly convex, the corresponding polyhedral partition is unique, providing a continuous parametric solution across the parameter space. However, implementing this algorithm for primal degenerate problems is not straightforward: a highly degenerate primal problem results in a high-dimensional optimal solution set for its dual. This is especially problematic for QPs (the auxiliary objective function) because CRs must be constructed from both primal and dual feasibility conditions, which amounts to projecting an H-representation (see \citet{Bemporad2001} for an overview of polyhedral representations) of a high-dimensional polyhedron onto a low-dimensional parameter space \cite{Tondel2003}. This is an NP-hard problem \cite{Tiwary2008}, so the algorithm of \citet{Spjotvold2007} does not scale well with problem size for highly-degenerate primal problems. Note that this issue is only relevant for QPs (and nonlinear programs in general) because dual feasibility is not implied by primal feasibility \cite{Nocedal2006}. In contrast, when a finite optimal primal solution exists, the LP strong duality theorem \cite{Sierksma2001} ensures that the dual problem always has an optimal feasible solution; thus, CRs can be constructed from the primal feasibility conditions alone without needing the forgoing polyhedral projection.  

Using a hierarchy of linear objective functions (\ie, a lexicographic LP) is an alternative technique to resolve optimal-solution multiplicity. This technique particularly suits dynamic flux balance analysis (DFBA), for example, where only the continuity of a subset of decision variables is of interest. DFBA involves dynamical systems with an LP embedded; the state variables depend on the optimal solution of a few decision variables in the LP, and the LP is parametrized with a few state variables \cite{Gomez2014, Harwood2015}. \citet{Gomez2016} applied this approach to study the dynamics of multi-species microbial consortia. Auxiliary objective functions were chosen to ensure a continuous solution for state variables corresponding to external metabolites that affect the extracellular dynamics. Alternatively, explicit solutions of the embedded LPs can be constructed using an MPP algorithm with a hierarchy of linear objectives, and, thus, the complexities of polyhedral projection can be avoided. 

Primal degeneracy, which is related but not equivalent to dual multiplicity, is the second common difficulty for MMP algorithms. Non-degeneracy implies that a balanced system of equations arises from the active set, from which a unique optimal solution can be readily constructed. Conversely, degeneracy implies that the active set corresponds to an overdetermined (linearly dependent) system of equations, and there are redundant active constraints\footnote{Not to be confused with redundant equality constraints with respect to a polyhedral feasible set (see Definition 2.1 of \citet{Telgen1982}). Unlike redundant equality constraints, eliminating a redundant active constraint can change the feasible set. A redundant active constraint is defined for an optimal solution and can become non-redundant if parameters are perturbed.}. To handle primal degeneracy, it is usually recommended to derive a balanced system from the active set by eliminating the redundant constraints \cite{Bemporad2002, Borrelli2003, Spjotvold2007}. However, this is algorithmically undesirable since the algorithm complexity grows combinatorially with problem size.

\citet{Jones2007} introduced a technique to handle degeneracy and multiplicity using pivoting and basis-updating rules from the Simplex algorithm through lexicographic perturbations of primal and dual constraints. To resolve primal degeneracy, primal constraints are weighted according to a prescribed priority by unevenly perturbing their right-hand sides, such that the corresponding optimal solution of the perturbed problem is non-degenerate \cite{Mangasarian1979a, Fukuda1997}. Similarly, primal multiplicity is handled by removing dual degeneracy using disparate weights for dual constraints through an uneven perturbation of the cost-vector components. This technique furnishes a unique polyhedral partition and continuous parametric solution for all decision variables for a prescribed priority order of primal and dual constraints. However, it cannot be used for problems where CR facets are shared between more than two CRs \cite{Spjotvold2006}. The algorithm of \citet{Jones2007} was later extended \cite{Jones2006} to address this limitation. As will be discussed in Section~\ref{sec:multiple}, this technique can be regarded as a lexicographic LP, where multiplicity is handled by selecting a vertex of the optimal face that satisfies a prescribed priority order the objective functions impose on the decision variables. 

In this paper, we modify the algorithms of \citet{Bemporad2002} and \citet{Borrelli2003} to improve the handling of degeneracy and primal multiplicity for LPs. We adopt a standard-form formulation \cite{Bertsimas1997} where the decision variables are partitioned into zero and non-zero components. This furnishes a natural generalization of non-basic and basic variables, which are well-known in linear programming, by unifying the treatment of degenerate and non-degenerate problems. It is also consistent with active-set-based MPP algorithms since zero variables correspond to active inequality constraints in standard form, which must remain unchanged in each CR by definition. Thus, computations are restricted to non-zero variables in each CR, significantly reducing the computational costs for highly degenerate problems where the number of zero variables is much larger than non-basic variables (\eg, metabolic-network models). Moreover, we use generalized inverses \cite{Ben-Israel2003} to handle degeneracy. This allows the construction of optimal solutions for primal degenerate problems without needing to eliminate redundant active constraints. Moreover, generalized inverses have a desirable algebraic structure, conforming to LU-decomposition techniques for full-rank and rank-deficient sparse matrices; thus, they can be efficiently computed using fast algorithms from sparse linear algebra. To resolve primal multiplicity, we apply auxiliary linear objective functions, as discussed above, and compare the resulting polyhedral partitions with those obtained from the algorithm of \citet{Spjotvold2007}. 

\section{Notation} \label{sec:notation}
Lowercase and uppercase boldface letters denote vectors and matrices, respectively. $\mathbf{I}_n$ and $\mathbf{0}$ denote an $n\times n$ identity matrix and a zero matrix of appropriate size, respectively. The dimension of an $m\times n$ zero matrix is explicitly stated as $\mathbf{0}_{m\times n}$ to avoid confusion when necessary. The all-one matrices $\mathbf{1}$ and $\mathbf{1}_{m\times n}$ are similarly defined. Boldface letters with subscript denote sub-matrices and sub-vectors of the original  matrix or vector that is defined in the context. Given an index set $\mathcal{I}$, $\mathbf{M}_{\mathcal{I}}$ is a sub-matrix comprising rows of $\mathbf{M}$ with the corresponding indices in $\mathcal{I}$. A covariant representation is adopted for vectors and a mixed covariant-contravariant for matrices; thus, $M_i^j$ refers to the element at the $i$th row and $j$th column of a matrix $\mathbf{M}$, and $v_i$ refers to the $i$th component of a column vector $\mathbf{v}$. The Einstein summation convention is used for brevity, where repeating a dummy index as a subscript and superscript imply summation over an appropriate range of the index. Polyhedral sets are denoted by uppercase Greek or normal letters. $\partial P$ denotes the boundary of the set $P$, and $\partial P_{1,2}$ denotes the common facet between two neighboring sets $P_1$ and $P_2$.

\section{Overview of MPP Algorithms Based on Active-Set Methods} \label{sec:overview}
An informal overview of active-set-based MPP algorithms is presented in this section. Main concepts underlying the algorithms of \citet{Bemporad2002}, \citet{Tondel2003a}, and \citet{Spjotvold2005} for LPs and QPs are demonstrated through illustrative examples. Two notions of degeneracy arise in the forthcoming discussion (see Definition 12.4 and Section 16.4 of \citet{Nocedal2006}): (i) Violation of linear independence constraint qualification (LICQ) and (ii) violation of strict complementarity. Unless explicitly stated, degeneracy in this paper always refers to the first case. We first discuss the degeneracy-multiplicity relationship in LPs as it plays a central role in the construction of continuous parametric solutions.

\begin{figure} 
	\centering
	\includegraphics[width=\linewidth]{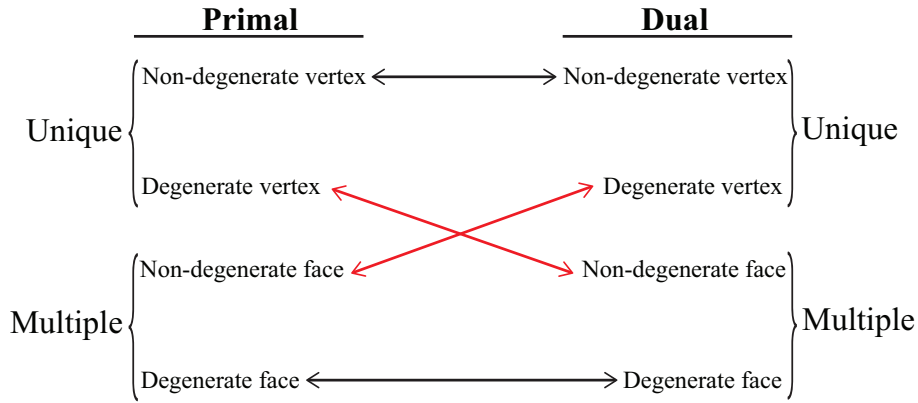}
  \caption{Schematic representation of Theorem 5 and Corollary 1 of \citet{Tijssen1998} concerning the relationship between multiplicity and degeneracy of LPs in canonical form. Double arrows indicate \emph{if and only if} relationships.}
	\label{fig:figure1}
\end{figure} 

\subsection{Degeneracy and Multiplicity in LPs} \label{sec:deg-mult}
Identifying the uniqueness or multiplicity of optimal solutions is a key step in MPP algorithms. Dual degeneracy is often mistakenly taken as being equivalent to primal multiplicity in the literature and, thereby, used as a criterion to determine whether an optimal basic feasible solution is unique or a vertex of an optimal face. \citet{Tijssen1998} provided a comprehensive account of degeneracy and multiplicity in LPs (schematically summarized in Fig.~\ref{fig:figure1}), showing that dual degeneracy of a vertex is only a necessary condition for primal multiplicity. They generalized the notion of a degenerate vertex to polyhedral faces by defining the degeneracy degree of a face $F$ of a polyhedron $P\subset\mathbbm{R}^n$ as
\begin{equation}%--------------------------------------------------------------------
    	\sigma(F,P)=\mathrm{bnd}(F,P)+\mathrm{dim}(F)-n,
	 \label{eqn:eq1}
\end{equation}%--------------------------------------------------------------------
where $\mathrm{bnd}(F,P)$ and $\mathrm{dim}(F)$ respectively denote the number of constraints in $P$ that are active at every point of $F$ and the dimension of $F$. Accordingly, $F$ is degenerate if $\sigma(F,P)>0$ and non-degenerate if $\sigma(F,P)=0$. Consequently, a face is degenerate if all its proper subsets are degenerate, and it is non-degenerate if it has at least one non-degenerate proper subset\footnote{Note that $F$ is not necessarily degenerate even if all its vertices are degenerate. To be a degenerate face, every (vertex and non-vertex) point of $F$ must be degenerate.}. Figure~\ref{fig:figure1} clearly demonstrates that, if a dual problem has a non-degenerate optimal face with degenerate and non-degenerate vertices, then the primal problem has a unique degenerate optimal solution. An LP solver may return a degenerate or non-degenerate dual vertex; but, primal multiplicity cannot be deduced in either case.  

\subsection{MPP Algorithms for QPs} \label{sec:mpp-qp}
Consider the quadratic programming problem
\begin{equation}%--------------------------------------------------------------------
		\left\{\begin{array}{l}
         	  z^{*}(\pmb{\theta})=\min\limits_{\mathbf{x}\in\mathbbm{R}^n} \frac{1}{2}\mathbf{x}^{\mathrm{T}}\mathbf{H}\mathbf{x}\\
						\mathbf{Gx}\le\mathbf{w}+\mathbf{F}\pmb{\theta}
		\end{array} \right.
	 \label{eqn:eq2}
\end{equation}%--------------------------------------------------------------------
with $\mathbf{x}\in\mathbbm{R}^n$, $\mathbf{H}\in\mathbbm{R}^{n\times n}$, $\mathbf{G}\in\mathbbm{R}^{m\times n}$, $\mathbf{w}\in\mathbbm{R}^m$, $\mathbf{F}\in\mathbbm{R}^{m\times q}$, $\pmb{\theta}\in\mathbbm{R}^q$, and 
\begin{equation}%--------------------------------------------------------------------
		\mathbf{H}=\left[\begin{array}{c c}
         	  1 & 0\\
						0 & 1
		\end{array} \right], \mathbf{G}=\left[\begin{array}{c c}
         	  0 & 1\\
						1 & 0\\
					 -1 & 0\\
					  0 & -1\\
					 -1 & -1
		\end{array} \right], \mathbf{w}=\left[\begin{array}{c}
         	  2\\
						3\\
					 -1\\
					 -1\\
					 -2
		\end{array} \right], \mathbf{F}=\left[\begin{array}{c c}
         	  0 & 0\\
						0 & 0\\
					  0 & 0\\
					 -1 & 0\\
					  0 & -1
		\end{array} \right],
	 \label{eqn:eq3}
\end{equation}%--------------------------------------------------------------------
where we seek a parametric solution pair $(\mathbf{x}^{*},\pmb{\mu}^{*})(\pmb{\theta})$ for $\pmb{\theta}\in[0,1]^2$. Here, $\pmb{\mu}$ denotes dual variables. For the given matrices in Eq.~(\ref{eqn:eq3}), $n=2$, $m=5$, and $q=2$. Moreover, provided LICQ is satisfied, the solution pair $(\mathbf{x}^{*},\pmb{\mu}^{*})$ must be a continuous piecewise affine function over a convex polyhedral partition of the parameter space (Theorem 2 of \citet{Bemporad2002}) because the Hessian matrix $\mathbf{H}$ is positive definite and the constraints are affine in $\pmb{\theta}$; when LICQ is violated, $\mathbf{x}^{*}$ still remains unique, but $\pmb{\mu}^{*}$ may have multiple solutions (see Section 12.3 of \citet{Nocedal2006}).

\begin{figure} 
	\centering
	\includegraphics[width=\linewidth]{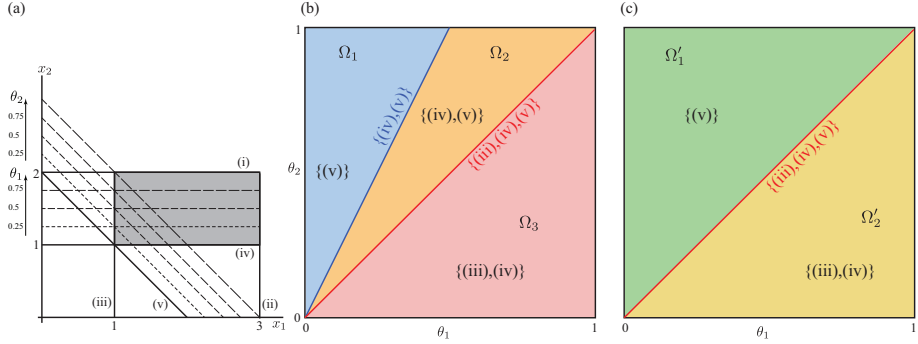}
  \caption{Schematic representation of the quadratic and linear programming problems given by Eqs.~(\ref{eqn:eq2}), (\ref{eqn:eq3}), and (\ref{eqn:eq8}) : (a) The base feasible set (shaded area) at $\pmb{\theta}=\mathbf{0}$ with the polyhedral partition of the parameter space for the (b) quadratic and (c) linear problem. Numbers enclosed by curly brackets denote the optimal active set in each region.}
	\label{fig:figure2}
\end{figure}

Polyhedral partitions of the parameter space (Fig.~\ref{fig:figure2}b) for the problem given by Eqs.~(\ref{eqn:eq2}) and (\ref{eqn:eq3}) (Fig.~\ref{fig:figure2}a) can be constructed using active-set-based algorithms \cite{Bemporad2002, Tondel2003a}, which are based on the Karush-Kuhn-Tucker (KKT) conditions
\begin{align}%--------------------------------------------------------------------
	\mathbf{H}\mathbf{x}+\mathbf{G}^{\mathrm{T}}\pmb{\mu}&=\mathbf{0}, \label{eqn:eq4}\\
	\mu_i(G_i^j x_j-w_i-F_i^j\theta_j)&=0, \quad i=1,\cdots,m\label{eqn:eq5}\\
	\mathbf{Gx}-\mathbf{w}-\mathbf{F}\pmb{\theta}&\le\mathbf{0},	\label{eqn:eq6}\\
	\pmb{\mu}&\ge\mathbf{0}.	\label{eqn:eq7}
\end{align}%--------------------------------------------------------------------
Equations~(\ref{eqn:eq4})-(\ref{eqn:eq7}) furnish the necessary and sufficient conditions to identify optimal solution pairs $(\mathbf{x}^{*},\pmb{\mu}^{*})$ for strictly convex QPs \cite{Nocedal2006}, such as
 Eqs.~(\ref{eqn:eq2}) and (\ref{eqn:eq3}). These algorithms recursively apply KKT conditions to find all the CRs covering a given search region $R^{(1)}$ in the parameter space ($R^{(1)}=[0,1]^2$ in Fig.~\ref{fig:figure2}b). Hereafter, we denote the search region in cycle $k$ by $R^{(k)}$ with $R^{(1)}$ the initial search region in which the CRs are to be constructed. Each cycle $k$ comprises five major steps: (i) Solve Eqs.~(\ref{eqn:eq4})-(\ref{eqn:eq7}) for $(\mathbf{x}^{*},\pmb{\mu}^{*})$ at $\pmb{\theta}_0$, which is usually chosen to be the Chebyshev center of $R^{(k)}$ \cite{Bemporad2002}; (ii) Determine the active set $\mathcal{A}(\mathbf{x}^{*}(\pmb{\theta}_0),\pmb{\theta}_0)$, inactive set $\mathscr{N}(\mathbf{x}^{*}(\pmb{\theta}_0),\pmb{\theta}_0)$, and optimal active set $\mathcal{A}^{*}(\pmb{\theta}_0):=\bigcap_{\mathbf{x}^{*}\in X^{*}(\pmb{\theta}_0)}\mathcal{A}(\mathbf{x}^{*},\pmb{\theta}_0)$, where $\mathcal{A}(\mathbf{x},\pmb{\theta}):=\{i|G_i^j x_j-w_i-F_i^j\theta_j=0\}$, $\mathscr{N}(\mathbf{x},\pmb{\theta}):=\{i|G_i^j x_j-w_i-F_i^j\theta_j<0\}$, and $X^{*}(\pmb{\theta})$ denotes the set of optimal primal solutions (optimal face) for any $\pmb{\theta}\in R^{(1)}$; (iii) Derive a balanced system of equations from active constraints $\mathbf{G}_{\mathcal{A}(\mathbf{x}^{*}(\pmb{\theta}_0),\pmb{\theta}_0)}$ and construct a parametric solution $\mathbf{x}^{*}(\pmb{\theta})$; (iv) Obtain an H-representation of the CR in the current cycle $\Omega^{(k)}$ from primal feasibility using inactive constraints $\mathbf{G}_{\mathscr{N}(\mathbf{x}^{*}(\pmb{\theta}_0),\pmb{\theta}_0)}$ and dual feasibility Eq.~(\ref{eqn:eq7}); (v) Search the remainder of the search region $\bar{R}^{(k)}:=R^{(k)}\setminus \Omega^{(k)}$ to find other CRs according to a search strategy. 

Two major search strategies were proposed by \citet{Bemporad2002} (Strategy  I) and \citet{Tondel2003a} (Strategy  II). In Strategy I, given an $f$-sided $\Omega^{(k)}$ and a search region $R^{(k)}$, $\bar{R}^{(k)}$ is partitioned into $f$ sub-regions, which are new regions to be consecutively searched for additional CRs. This procedure results in a search tree, any branch of which terminates when a search region $R^{(l)}$ with $\bar{R}^{(l)}=\emptyset$ is encountered. However, partitioning $\bar{R}^{(k)}$ can split a CR into several sub-regions, each requiring a separate calculation of the active set and optimal solutions. Therefore, one may need to execute Steps (i)-(v) several times more than there are CRs. In Strategy II, assuming that common facets are not shared by more than two CRs, the active set in each neighboring CR of $\Omega^{(k)}$ is determined from the active set in CR and the type of their common facets. This strategy avoids splitting CRs, and, thus, Steps (i)-(v) are executed as many times as there are CRs. Here, facet types are defined either by a dual-feasibility or primal-feasibility constraint, and LICQ can only be violated on the facets \cite{Tondel2003a}. Hence, this strategy cannot be applied to many practical problems, such as FBA, that are highly degenerate in large regions of their parameter space (see Chapter 3 of \citet{Spjoetvold2008} for an overview of various search strategies). 

\begin{table}[t]
\caption{Primal solution $\mathbf{x}^{*}$, dual solution $\pmb{\mu}^{*}$, and optimal active set $\mathcal{A}^{*}$ in critical regions and on critical-region boundaries of Fig.~\ref{fig:figure2}b.}
\centering
\begin{tabular}{c @{\hskip.1cm} c @{\hskip.1cm} c @{\hskip.1cm} c @{\hskip.05cm} c}
\hline
\shortstack[c]{Region} & \shortstack[c]{$\mathcal{A}^{*}$}  & \shortstack[c]{$\mathbf{x}^{*}$} & \shortstack[c]{$\pmb{\mu}^{*}$} & \shortstack[c]{Feasible set} \\
\hline
$\Omega_1$ & $\{\rm{(v)}\}$ & $\left[\arraycolsep=.2pt\begin{array}{c}
         	  1+\theta_2/2\\
						1+\theta_2/2
		\end{array} \right]$ & {\arraycolsep=.1pt$
						\left[\begin{array}{c}
						\mathbf{0}_{4\times1}\\
						1+\theta_2/2
						\end{array}\right]$} & Fig.~\ref{fig:figure3}e\\
$\Omega_2$ & $\{\rm{(iv)},\rm{(v)}\}$ & $\left[\arraycolsep=.2pt\begin{array}{c}
         	  1+\theta_2-\theta_1\\
						1+\theta_1
		\end{array} \right]$ & {\arraycolsep=.1pt$
						\left[\begin{array}{c}
						\mathbf{0}_{3\times1}\\
						2\theta_1-\theta_2\\
						1+\theta_2-\theta_1
						\end{array}\right]$} & Fig.~\ref{fig:figure3}d\\
$\Omega_3$ & $\{\rm{(iii)},\rm{(iv)}\}$ & $\left[\arraycolsep=.2pt\begin{array}{c}
         	  1\\
						1+\theta_1
		\end{array} \right]$ & {\arraycolsep=.1pt$
						\left[\begin{array}{c}
						\mathbf{0}_{2\times1}\\
						1\\
						1+\theta_1\\
						0
						\end{array}\right]$} & Fig.~\ref{fig:figure3}c\\
$\partial \Omega_{2,3}$ & $\{\rm{(iii)},\rm{(iv)},\rm{(v)}\}$ & $\left[\arraycolsep=.2pt\begin{array}{c}
         	  1\\
						1+\theta_1
		\end{array} \right]$ & {\arraycolsep=.1pt$
						\left[\begin{array}{c}
						\mathbf{0}_{2\times1}\\
						(1-\theta_1)/3\\
						(1+2\theta_1)/3\\
						(2+\theta_1)/3
						\end{array}\right]+\left[\begin{array}{c}
						\mathbf{0}_{2\times1}\\
						1\\
						1\\
						-1
						\end{array}\right]\mu_z$} & Fig.~\ref{fig:figure3}b\\
$\partial \Omega_{1,2}$ & $\{\rm{(iv)},\rm{(v)}\}$ & $\left[\arraycolsep=.2pt\begin{array}{c}
         	  1+\theta_1\\
						1+\theta_1
		\end{array} \right]$ & {\arraycolsep=.1pt$
						\left[\begin{array}{c}
						\mathbf{0}_{4\times1}\\
						1+\theta_1
						\end{array}\right]$} & Fig.~\ref{fig:figure3}a\\

\hline
\end{tabular}
\label{tbl:table1}
\end{table}

\begin{figure} [t]
	\centering
	\includegraphics[width=\linewidth]{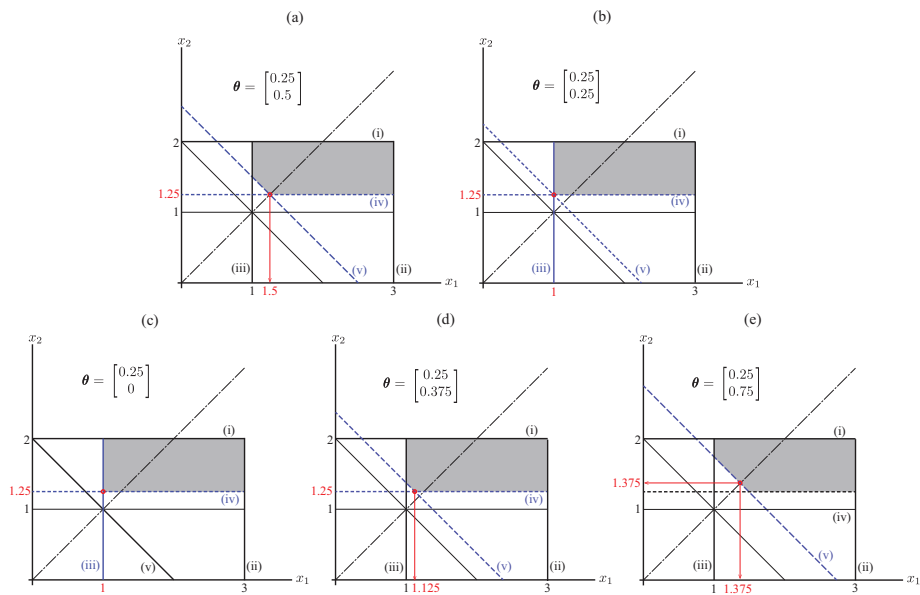}
  \caption{Changes in the feasible set (shaded area) and optimal primal solution $\mathbf{x}^{*}$ (open circle) of Fig.~\ref{fig:figure2} at fixed $\theta_1=0.25$ with (a) $\theta_2=0.5$, (b) $\theta_2=0.25$, (c) $\theta_2=0$, (d) $\theta_2=0.375$, and (e) $\theta_2=0.75$. The line $x_2=x_1$ (dash-dotted) is shown to help identify the optimal solution.}
	\label{fig:figure3}
\end{figure}

Applying MPP algorithms to the QP given by Eqs.~(\ref{eqn:eq2}) and (\ref{eqn:eq3}) using either search strategy yields three $q$-dimensional CRs and two $(q-1)$-dimensional facets\footnote{For a problem with $q$ parameters, a $q$-dimensional (full-dimensional) critical region is a polyhedron that cannot be embedded in an affine space with fewer dimensions than $q$. Unless stated otherwise, critical region always refers to a full-dimensional polyhedron in this paper.}, as shown in Fig.~\ref{fig:figure2}b. Table~\ref{tbl:table1} summarizes their respective parametric primal and dual solutions. Note that, in this example, the set of active constraints for each solution is also the optimal active set because the primal solution is unique in the entire parameter space. The problem is primal non-degenerate in all three critical regions. The primal solution is a non-vertex point in $\Omega_1$ (Fig.~\ref{fig:figure3}e) and a vertex point in $\Omega_2$ and $\Omega_3$ (Fig.~\ref{fig:figure3}d and c). Particular attention should be paid to the solution form and active set on $\partial \Omega_{1,2}$ and $\partial \Omega_{2,3}$. Hyperplanes containing the common facets (and any lower-dimensional common faces) of neighboring critical regions are loci of points at which active-set exchanges occur, and consequently, global parametric solutions of LPs and QPs could become non-smooth or discontinuous. Along any path from $\Omega_2$ to $\Omega_3$, the constraints (iii) and (v) respectively enter and leave the active set at a point on $\partial \Omega_{2,3}$. In contrast, along any path from $\Omega_1$ to $\Omega_2$, only the constraint (v) enters the active set. The primal solution is degenerate on $\partial \Omega_{2,3}$ (Fig.~\ref{fig:figure3}b) and non-degenerate on $\partial \Omega_{1,2}$ (Fig.~\ref{fig:figure3}a). Therefore, the dual solution is non-unique  (one dimensional here) in the former and unique in the latter case (see Table~\ref{tbl:table1}).  

The polyhedral partition in Fig.~\ref{fig:figure2}b features the two types of hyperplanes separating neighboring critical regions, as discussed by \citet{Tondel2003a}, namely, those corresponding to (i) primal feasibility (Type I, resulting from Eq.~(\ref{eqn:eq6})) and (ii) dual feasibility (Type II, resulting from Eq.~(\ref{eqn:eq7})) constraints. The facet $\partial \Omega_{2,3}$ is a hyperplane of Type I with degenerate primal solutions, separating two full-dimensional critical regions with non-degenerate primal solutions. Crossing this facet from either side in the parameter space, the vertex corresponding to the active set in one CR becomes infeasible in the other. For example, the vertex at the intersection of the constraints (iv) and (v) is optimal in $\Omega_2$ (Fig.~\ref{fig:figure3}d), but infeasible in $\Omega_3$ (Fig.~\ref{fig:figure3}c). Note that the degeneracy of the primal solution in neighboring CRs and their separating facets can generally be of any higher degree. Therefore, the consequences of Theorem 2 and Lemma 1 of \citet{Tondel2003a}, which are central to the search Strategy II discussed above, are not always applicable. In contrast, the facet $\partial \Omega_{1,2}$ is a hyperplane of Type II with non-degenerate primal solutions, separating two full-dimensional critical regions with non-degenerate primal solutions. Crossing $\partial \Omega_{1,2}$ from $\Omega_2$ to $\Omega_1$, the optimal vertex at the intersection of the constraints (iv) and (v) remains primal feasible (Fig.~\ref{fig:figure3}e), but is no longer a dual feasible solution. Along any such path, the constraint (iv) first becomes weakly active ($\mu^{*}_4=0$ and $\rm{(iv)}\in\mathcal{A}$) on $\partial \Omega_{1,2}$ and then inactive in $\Omega_1$.

\begin{figure} [t]
	\centering
	\includegraphics[width=\linewidth]{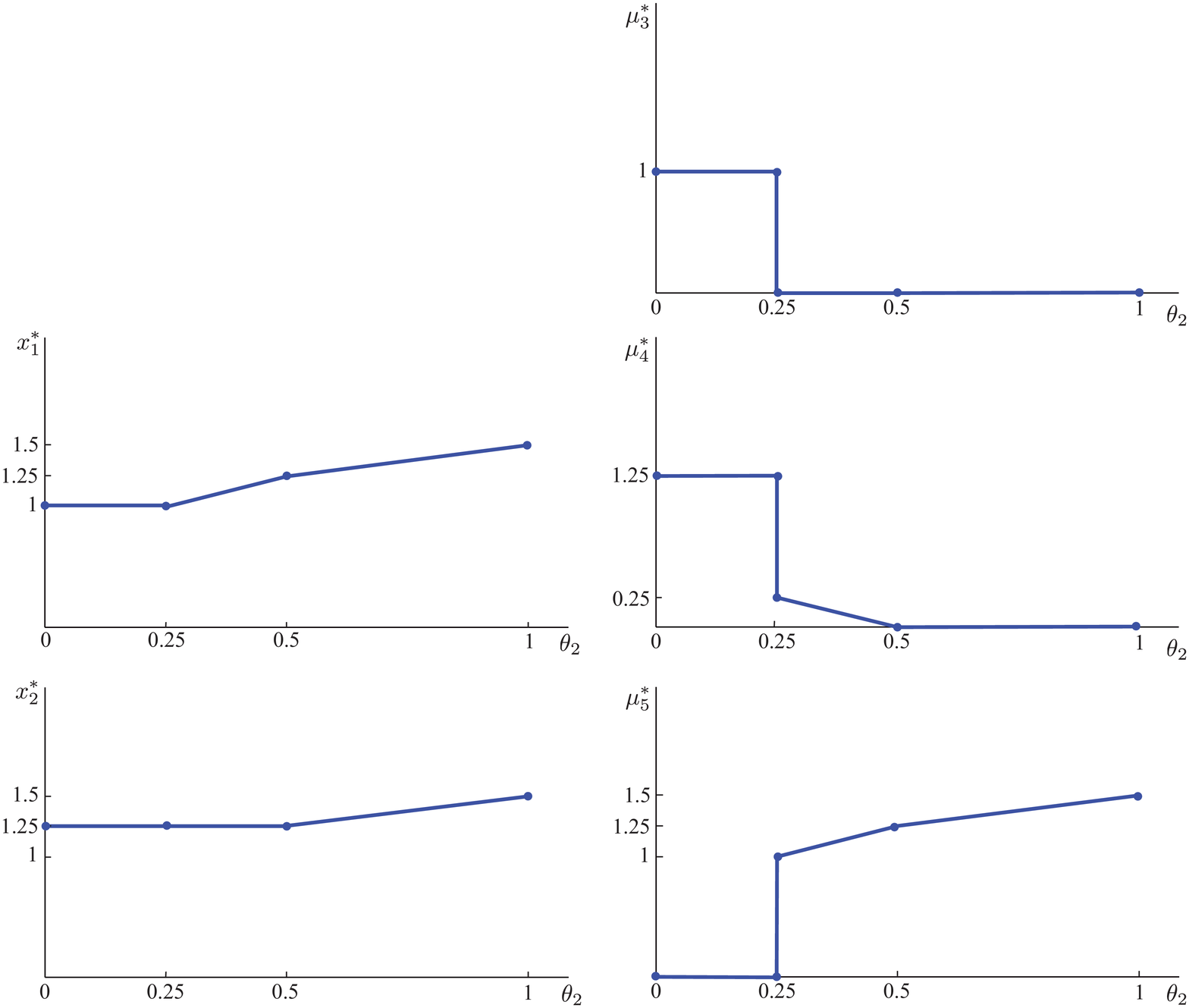}
  \caption{Continuity of optimal primal $\mathbf{x}^{*}$ and dual $\pmb{\mu}^{*}$ solutions of Fig.~\ref{fig:figure2} along the path $\mathcal{C}=\{\pmb{\theta}\in[0,1]^2|\theta_1=0.25\}$.}
	\label{fig:figure4}
\end{figure}

Primal degeneracy also has important implications for the continuity of parametric solutions. Figure~\ref{fig:figure4} provides a qualitative comparison of optimal solutions along a parameter-space path crossing degenerate and non-degenerate facets. Here, the optimal primal and dual solutions are plotted along a vertical path ($\theta_1=0.25$) in Fig.~\ref{fig:figure2}b, which crosses the degenerate facet $\partial \Omega_{2,3}$ at $\theta_2=0.25$ and non-degenerate facet $\partial \Omega_{1,2}$ at $\theta_2=0.5$. Here, the positive definite $\textbf{H}$ guarantees a continuous $\textbf{x}^{*}(\pmb{\theta})$ over the entire parameter space (Fig.~\ref{fig:figure4}, left panel). However, $\pmb{\mu}^{*}(\pmb{\theta})$ is non-unique on $\partial \Omega_{2,3}$ due to primal degeneracy (Fig.~\ref{fig:figure4}, right panel) and, therefore, is discontinuous with respect to the Hausdorff metric \cite{Mangasarian1987}. Note that the solution set of $\pmb{\mu}^{*}$ on $\partial \Omega_{2,3}$ is bounded since the Mangasarian-Fromovitz constraint qualification holds in this example (see Section 12.6 of \citet{Nocedal2006}). As will be discussed in Section~\ref{sec:improvement}, compared to LPs, dual multiplicity poses major computational challenges to MPP algorithms for QPs. 

Next, consider the linear programming problem
\begin{equation}%--------------------------------------------------------------------
		\left\{\begin{array}{l}
         	  z^{*}(\pmb{\theta})=\min\limits_{\mathbf{x}\in\mathbbm{R}^n} \mathbf{c}^{\mathrm{T}}\mathbf{x}\\
						\mathbf{Gx}\le\mathbf{w}+\mathbf{F}\pmb{\theta}
		\end{array} \right.,
	 \label{eqn:eq8}
\end{equation}%--------------------------------------------------------------------
with the same feasible set as given by Eq.~(\ref{eqn:eq3}) and $\mathbf{c}=[1 \quad 1]^{\mathrm{T}}$. Note that LPs can be viewed as QPs with positive-semidefinite Hessians. Thus, the uniqueness of the primal solution is no longer guaranteed. Figure~\ref{fig:figure2}c shows the polyhedral partition arising from this LP. The primal solution is unique in $\Omega'_2$ and on $\partial \Omega'_{1,2}$, which are identical to those of $\Omega_3$ and $\partial \Omega_{2,3}$ in Fig.~\ref{fig:figure2}b. However, there are multiple primal solutions in $\Omega'_1$, all lying on a face of the feasible set corresponding to the constraint (v). Therefore, the optimal active set in $\Omega'_1$ is $\mathcal{A}^{*}=\{\rm{(v)}\}$. Nevertheless, to construct an explicit parametric solution, one needs to choose a point from the solution set according to a well-defined rule. For example, \citet{Borrelli2003} suggested to arbitrarily choose one vertex of the optimal face and construct the respective CR. Following this rule, either $\mathcal{A}=\{\rm{(iv)},\rm{(v)}\}$ or $\mathcal{A}=\{\rm{(iii)},\rm{(v)}\}$ can be selected as the active set corresponding to the chosen optimal vertex in $\Omega'_1$, both furnishing a continuous parametric solution in the entire parameter space. However, this strategy does not always lead to continuous parametric solutions (see \citet{Spjotvold2005} for an example).

Using auxiliary objective functions is another technique to construct unique parametric solutions systematically. \citet{Spjotvold2007} proposed to apply an auxiliary norm-minimizing objective function when primal multiplicity is encountered in LPs. The continuity of the parametric primal solution is ensured by the positive definiteness of the Hessian. However, this method inherits the computational difficulties of parametric QPs in handling  dual multiplicity. Applying this technique to Eq.~(\ref{eqn:eq8}), the parameter space is partitioned into three CRs that are identical to those of Eq.~(\ref{eqn:eq2}) shown in Fig.~\ref{fig:figure2}b. Here, the auxiliary objective function splits $\Omega'_1$ into two sub-regions ($\Omega_1$ and $\Omega_2$ in Fig.~\ref{fig:figure2}b), leaving $\Omega'_2$ unaffected (identical to $\Omega_3$ in Fig.~\ref{fig:figure2}b). Alternatively, a hierarchy of linear auxiliary objective functions can be used to resolve primal multiplicity. For example, an auxiliary linear objective with $\mathbf{d}_1=[-1 \quad 1]^{\mathrm{T}}$ or $\mathbf{d}_1=[1 \: -1]^{\mathrm{T}}$ leads to a unique parametric solution in the entire parameter space for Eq.~(\ref{eqn:eq8}). These, respectively, yield $\mathcal{A}=\{\rm{(iii)},\rm{(v)}\}$ or $\mathcal{A}=\{\rm{(iv)},\rm{(v)}\}$ in $\Omega'_1$ without changing the polyhedral partition of the parameter space (the same as Fig.~\ref{fig:figure2}c).

\section{Improvements on Active-Set-Based MPP Algorithms} \label{sec:improvement}
In this section, we provide improvements to existing MPP methods \cite{Bemporad2002, Borrelli2003, Spjotvold2007} for LPs, specifically focusing on the Steps (iii) and (iv) of active-set-based algorithms discussed in Section~\ref{sec:mpp-qp}. We adopt the search Strategy I for Step (v) since it does not rely on the degree of degeneracy of optimal solutions in CRs and their facets. Emphasis is placed on the handling of degeneracy and multiplicity as they frequently arise in practical applications, such as FBA. We adopt a standard-form formulation and consider the problem 
\begin{equation}%--------------------------------------------------------------------
		\left\{\begin{array}{l}
         	  z^{*}(\pmb{\theta})=\max\limits_{\mathbf{x}\in\mathbbm{R}^n} \mathbf{c}^{\mathrm{T}}\mathbf{x}\\
						\mathbf{Ax}=\mathbf{w}+\mathbf{F}\pmb{\theta}\\
						\mathbf{x}\geq\mathbf{0}
		\end{array} \right.
	 \label{eqn:eq9}
\end{equation}%--------------------------------------------------------------------
with $\mathbf{x}\in\mathbbm{R}^n$, $\mathbf{c}\in\mathbbm{R}^{n}$, $\mathbf{A}\in\mathbbm{R}^{m\times n}$, $\mathbf{w}\in\mathbbm{R}^m$, $\mathbf{F}\in\mathbbm{R}^{m\times q}$, $\pmb{\theta}\in R^{(1)}\subset\mathbbm{R}^q$, where the primary objective function is linear. We assume that the problem is well-conditioned, so that $[\mathbf{A}^{\mathrm{T}}\:\mathbf{c}]$ has full column rank. Recall, to construct a CR in each cycle of these algorithms, one needs a primal solution at a point of the search region $\mathbf{x}^{*}(\pmb{\theta}_0)$ to determine the optimal active set $\mathcal{A}^{*}(\pmb{\theta}_0)$ and active set $\mathcal{A}(\mathbf{x}^{*}(\pmb{\theta}_0),\pmb{\theta}_0)$. Therefore, in the remainder of this section, we assume $\mathbf{x}^{*}(\pmb{\theta}_0)$ is known for a point in the current search region $\pmb{\theta}_0\in R^{(k)}$, demonstrating how to determine a parametric solution of the primal and dual problems and construct the H-representation of the respective critical region $\Omega^{(k)}$. When the primal solution is non-unique, $\mathbf{x}^{*}(\pmb{\theta}_0)$ is assumed to be a relative interior point (\eg, Chebyshev center) of the optimal face, the construction of which is explained later in this section. 

The standard-form formulation has the advantage that non-negativities are the only inequality constraints so that active and inactive inequalities are associated with zero and non-zero decision variables, respectively. Moreover, the optimal active set $\mathcal{A}^{*}(\pmb{\theta})$ does not change $\forall \pmb{\theta}\in \Omega^{(k)}$ by definition, so zero (non-zero) variables must remain zero (non-zero) in the current CR. As a result, a generalization of the concept of basic and non-basic variables from linear programming naturally follows: decision variables are partitioned into non-zero and zero variables rather than basic and non-basic because it significantly simplifies the KKT system, unifying the treatment of degenerate and non-degenerate problems.   

Since equality constraints are always active, only the index set of the inequalities $\mathcal{J}:=\{1,2,\cdots,n\}$ needs to be introduced into the KKT formulation. Let $\{\mathcal{J}_1,\mathcal{J}_0\}$ be a partition of $\mathcal{J}$ with $\mathcal{J}_1:=\{j\in\mathcal{J}|\exists\mathbf{x}\in X^{*}(\pmb{\theta}_0), x_{j}(\pmb{\theta}_0)>0\}$, $\mathcal{J}_0:=\{j\in\mathcal{J}|\forall\mathbf{x}\in X^{*}(\pmb{\theta}_0), x_{j}(\pmb{\theta}_0)=0\}$, $n_1:=|\mathcal{J}_1|$, and $n_0:=|\mathcal{J}_0|$. The respective partition of the decision variables is $(\mathbf{x}_{\mathcal{J}_1},\mathbf{x}_{\mathcal{J}_0})$. A permutation matrix $\mathbf{P}$ can be defined from $\mathbf{x}^{*}(\pmb{\theta}_0)$, which we use to transform the current coordinate system into one with respect to which all zero and non-zero variables and their respective columns in the constraint matrix $\mathbf{A}$ are completely separated (see \ref{sec:appendixa}). The index map $\xi$ associated with $\mathbf{P}$ provides the index set $\bar{\mathcal{J}}=\xi(\mathcal{J})$ and its respective partition $\{\bar{\mathcal{J}}_1,\bar{\mathcal{J}}_0\}$, decision-variable partition $(\bar{\mathbf{x}}_{1},\bar{\mathbf{x}}_{0})$, cost-vector partition $(\bar{\mathbf{c}}_{1},\bar{\mathbf{c}}_{0})$, and dual-variable partition $(\bar{\pmb{\mu}}_{1},\bar{\pmb{\mu}}_{0})$ in the transformed coordinate system with $\bar{\mathbf{x}}_{1}:=\bar{\mathbf{x}}_{\bar{\mathcal{J}}_1}=\mathbf{x}_{\Xi(\bar{\mathcal{J}}_1)}$, $\bar{\mathbf{x}}_{0}:=\bar{\mathbf{x}}_{\bar{\mathcal{J}}_0}=\mathbf{x}_{\Xi(\bar{\mathcal{J}}_0)}$, $\bar{\mathbf{c}}_{1}:=\bar{\mathbf{c}}_{\bar{\mathcal{J}}_1}=\mathbf{c}_{\Xi(\bar{\mathcal{J}}_1)}$, $\bar{\mathbf{c}}_{0}:=\bar{\mathbf{c}}_{\bar{\mathcal{J}}_0}=\mathbf{c}_{\Xi(\bar{\mathcal{J}}_0)}$, $\bar{\pmb{\mu}}_{1}:=\bar{\pmb{\mu}}_{\bar{\mathcal{J}}_1}$, and $\bar{\pmb{\mu}}_{0}:=\bar{\pmb{\mu}}_{\bar{\mathcal{J}}_0}$. The transformed constraint matrix can be written $\bar{\mathbf{A}}=[\bar{\mathbf{A}}_1|\bar{\mathbf{A}}_0]$, where $\bar{\mathbf{A}}_1$ and $\bar{\mathbf{A}}_0$ are sub-matrices containing all the columns of $\mathbf{A}$ that correspond to non-zero and zero decision variables, respectively. Accordingly, Eq.~(\ref{eqn:eq9}) transforms to 
\begin{equation}%--------------------------------------------------------------------
		\left\{\begin{array}{l}
         	  z^{*}(\pmb{\theta})=\min\limits_{\bar{\mathbf{x}}\in\mathbbm{R}^n} -\bar{\mathbf{c}}^{\mathrm{T}}\bar{\mathbf{x}}\\
						\bar{\mathbf{A}}\bar{\mathbf{x}}=\mathbf{w}+\mathbf{F}\pmb{\theta}\\
						\bar{\mathbf{x}}\geq\mathbf{0}
		\end{array} \right.
	 \label{eqn:eq10}
\end{equation}%--------------------------------------------------------------------
with the respective KKT system
\begin{align}%--------------------------------------------------------------------
	\bar{\mathbf{c}}_{1}+\bar{\mathbf{A}}^{\mathrm{T}}_{1}\bar{\pmb{\lambda}}&=\mathbf{0}_{n_1\times 1}, \label{eqn:eq11}\\
	\bar{\mathbf{c}}_{0}+\bar{\mathbf{A}}^{\mathrm{T}}_{0}\bar{\pmb{\lambda}}+\bar{\pmb{\mu}}_{0}&=\mathbf{0}_{n_0\times 1}, \label{eqn:eq12}\\
	\bar{\mathbf{A}}_1\bar{\mathbf{x}}_1&=\mathbf{w}+\mathbf{F}\pmb{\theta},	\label{eqn:eq13}\\
	\bar{\mathbf{x}}_1&\ge\mathbf{0}_{n_1\times 1},	\label{eqn:eq14}\\
	\bar{\pmb{\mu}}_0&\ge\mathbf{0}_{n_0\times 1},	\label{eqn:eq15}
\end{align}%--------------------------------------------------------------------
where $\bar{\pmb{\lambda}}\in\mathbbm{R}^m$ and $\bar{\pmb{\mu}}\in\mathbbm{R}^n$ are the dual variables corresponding to equality and non-negativity constraints of Eq.~(\ref{eqn:eq10}). A solution triple $(\bar{\mathbf{x}}^{*}_1,\bar{\pmb{\lambda}}^{*},\bar{\pmb{\mu}}^{*}_0)$ satisfying Eqs.~(\ref{eqn:eq11})-(\ref{eqn:eq15}) is called a KKT point of Eq.~(\ref{eqn:eq10}) in $\Omega^{(k)}$. Note that $\bar{\mathbf{x}}^{*}_0=\mathbf{0}$ and $\bar{\pmb{\mu}}^{*}_1=\mathbf{0}$ identically hold $\forall \pmb{\theta}\in \Omega^{(k)}$, respectively, by definition and by complementary slackness; thus, they are eliminated from the forgoing KKT system. Consequently, complementary slackness is automatically satisfied by any solution triple of Eqs.~(\ref{eqn:eq11})-(\ref{eqn:eq15}), so it needs not be explicitly included in the KKT system. 

Parametric solutions in each CR can be characterized from two independent standpoints, namely, primal degeneracy and primal multiplicity. The combination of these leads to four possible cases, as highlighted by \citet{Spjotvold2007}, each requiring specific techniques to be handled by MPP algorithms. Noting that the solution characteristics do not change in each CR by definition, degeneracy and multiplicity criteria need to be checked only at $\mathbf{x}^{*}(\pmb{\theta}_0)$ for the current critical region $\Omega^{(k)}$. Primal degeneracy can readily be identified from the number of non-zero decision variables $n_1$. However, identifying primal multiplicity is not as straightforward. Dual degeneracy of a vertex is often mistakenly used in the literature as a criterion to determine whether a primal solution is unique; as discussed in Section~\ref{sec:deg-mult}, it is only a necessary condition for primal multiplicity (see Fig.~\ref{fig:figure1}). \citet{Appa2002} provided necessary and sufficient conditions for primal multiplicity, which are algorithmically convenient to implement, only requiring the solution of an auxiliary LP. Once the solution characteristics at $\mathbf{x}^{*}(\pmb{\theta}_0)$ have been identified, the following sections detail how to construct $(\bar{\mathbf{x}}^{*},\bar{\pmb{\lambda}}^{*},\bar{\pmb{\mu}}^{*})(\pmb{\theta})$ and the H-representation of $\Omega^{(k)}$ for each case.  

\subsection{Unique Primal Solution} \label{sec:unique}
\subsubsection{Non-Degenerate Case} \label{sec:unique-nd}
This case applies when $X^{*}(\pmb{\theta}_0)$ only contains a vertex of the feasible set and $n_1=m$. The primal and dual solutions can readily be constructed from the KKT conditions
\begin{align}%--------------------------------------------------------------------
	\bar{\mathbf{x}}^{*}_1&=\bar{\mathbf{A}}^{-1}_{1}(\mathbf{w}+\mathbf{F}\pmb{\theta}),	\label{eqn:eq16}\\
	\bar{\pmb{\lambda}}^{*}&=-(\bar{\mathbf{A}}^{\mathrm{T}}_{1})^{-1}\bar{\mathbf{c}}_{1}, \label{eqn:eq17}\\
	\bar{\pmb{\mu}}^{*}_0&=-\bar{\mathbf{c}}_{0}+\bar{\mathbf{A}}^{\mathrm{T}}_{0}(\bar{\mathbf{A}}^{\mathrm{T}}_{1})^{-1}\bar{\mathbf{c}}_{1}. \label{eqn:eq18}
\end{align}%--------------------------------------------------------------------
As expected, the dual solutions $\bar{\pmb{\lambda}}^{*}$ and $\bar{\pmb{\mu}}^{*}$ are unique and independent of $\pmb{\theta}$. Here, dual feasibility Eq.~(\ref{eqn:eq15}) is always guaranteed to hold, irrespective of $\pmb{\theta}$, by the strong duality theorem (Theorem 2.2 of \citet{Sierksma2001}), so $\Omega^{(k)}$ is defined only by primal feasibility Eq.~(\ref{eqn:eq14})  
\begin{equation}%--------------------------------------------------------------------
	\Omega^{(k)}=\{\pmb{\theta}\in R^{(k)}|\bar{\mathbf{A}}^{-1}_{1}(\mathbf{w}+\mathbf{F}\pmb{\theta})\ge\mathbf{0}\}. 
	 \label{eqn:eq19}
\end{equation}%--------------------------------------------------------------------

\subsubsection{Degenerate Case} \label{sec:unique-d}
Similarly to the previous case, $X^{*}(\pmb{\theta}_0)$ only contains a vertex of the feasible set, but primal degeneracy implies $n_1<m$. Consequently, the KKT system provides an overdetermined system of equations for $\bar{\mathbf{x}}^{*}_1$ and an underdetermined system for $\bar{\pmb{\lambda}}^{*}$. Overdetermined systems resulting from primal degeneracy are handled in the existing MPP algorithms by constructing a balanced system from Eq.~(\ref{eqn:eq13}) either by eliminating linearly dependent equalities \cite{Bemporad2002} or applying a reduction technique \cite{Borrelli2003}. As will be discussed in Section~\ref{sec:complexity}, the computational complexity these extra steps add to that of the non-degenerate case scales as $O(m!)$ in the former and $O(m^3)$ in the latter. QR factorization is a standard reduction technique for handling degeneracy \cite{Tondel2003}. Application of QR factorization to Eq.~(\ref{eqn:eq13}) is equivalent to a linear transformation of the decision variables, so one needs to reformulate the KKT system in terms of transformed primal and dual variables. Here, we propose an alternative technique based on generalized inverses \cite{Ben-Israel2003} with the same computational complexity as QR factorization, which can be directly applied to Eqs.~(\ref{eqn:eq11})-(\ref{eqn:eq15}). The Moore-Penrose inverse particularly suits our purposes because it is unique and satisfies most properties of standard inverses \cite{Ben-Israel2003}. Two variations we frequently use are the left pseudoinverse for overdetermined and right pseudoinverse for underdetermined systems. Accordingly, the primal solution is constructed using the left pseudoinverse of $\bar{\mathbf{A}}_{1}$
\begin{equation}%--------------------------------------------------------------------
	\bar{\mathbf{x}}^{*}_1=(\bar{\mathbf{A}}^{\mathrm{T}}_{1}\bar{\mathbf{A}}_{1})^{-1}\bar{\mathbf{A}}^{\mathrm{T}}_{1}(\mathbf{w}+\mathbf{F}\pmb{\theta}).
	 \label{eqn:eq20}
\end{equation}%--------------------------------------------------------------------

Since the dual variable $\bar{\pmb{\lambda}}^{*}$ has multiple solutions, it is decomposed into a minimum-norm and a null-space component as $\bar{\pmb{\lambda}}^{*}=\mathbf{Y}\bar{\pmb{\lambda}}^{*}_{y}+\mathbf{Z}\bar{\pmb{\lambda}}^{*}_{z}$ (see Section 15.3 of \citet{Nocedal2006}), where the columns of $\mathbf{Z}$ and $\mathbf{Y}$ span the null space of $\bar{\mathbf{A}}^{\mathrm{T}}_{1}$ and its orthogonal complement, respectively. The minimum-norm component $\bar{\pmb{\lambda}}^{*}_p:=\mathbf{Y}\bar{\pmb{\lambda}}^{*}_{y}$, to be regarded as the unique particular component of $\bar{\pmb{\lambda}}^{*}$, can also be derived from Eq.~(\ref{eqn:eq11}) using the right pseudoinverse of $\bar{\mathbf{A}}^{\mathrm{T}}_{1}$ as $\bar{\pmb{\lambda}}^{*}_p=-\bar{\mathbf{A}}_{1}(\bar{\mathbf{A}}^{\mathrm{T}}_{1}\bar{\mathbf{A}}_{1})^{-1}\bar{\mathbf{c}}_{1}$. Moreover, the degeneracy degree of the primal problem $m-n_1$ determines the dimension of the homogeneous (null-space)  component $\bar{\pmb{\lambda}}^{*}_h:=\mathbf{Z}\bar{\pmb{\lambda}}^{*}_{z}$. The dual solutions can be written
\begin{align}%--------------------------------------------------------------------
	\bar{\pmb{\lambda}}^{*}&=-\bar{\mathbf{A}}_{1}(\bar{\mathbf{A}}^{\mathrm{T}}_{1}\bar{\mathbf{A}}_{1})^{-1}\bar{\mathbf{c}}_{1}+\mathbf{Z}\bar{\pmb{\lambda}}^{*}_{z}, \label{eqn:eq21}\\
	\bar{\pmb{\mu}}^{*}_0&=-\bar{\mathbf{c}}_{0}+\bar{\mathbf{A}}^{\mathrm{T}}_{0}\bar{\mathbf{A}}_{1}(\bar{\mathbf{A}}^{\mathrm{T}}_{1}\bar{\mathbf{A}}_{1})^{-1}\bar{\mathbf{c}}_{1}-\bar{\mathbf{A}}^{\mathrm{T}}_{0}\mathbf{Z}\bar{\pmb{\lambda}}^{*}_{z}. \label{eqn:eq22}
\end{align}%--------------------------------------------------------------------
To ensure that dual feasibility holds in the KKT system, it is sufficient to find at least a $\bar{\pmb{\lambda}}^{*}_{z}\in\mathbbm{R}^{(m-n_1)}$ such that $\bar{\pmb{\mu}}^{*}_0\ge\mathbf{0}$. As in the previous case, this is ensured by the strong duality theorem (see \ref{sec:appendixb} for further clarification). Consequently, $\Omega^{(k)}$ is defined only by primal feasibility  
\begin{equation}%--------------------------------------------------------------------
	\Omega^{(k)}=\{\pmb{\theta}\in R^{(k)}|(\bar{\mathbf{A}}^{\mathrm{T}}_{1}\bar{\mathbf{A}}_{1})^{-1}\bar{\mathbf{A}}^{\mathrm{T}}_{1}(\mathbf{w}+\mathbf{F}\pmb{\theta})\ge\mathbf{0}\}.
	 \label{eqn:eq23}
\end{equation}%--------------------------------------------------------------------

A few remarks are warranted concerning the computations of the primal and dual solutions in this section. First, the sparsity of $\bar{\mathbf{A}}_{1}$ for metabolic-network models can be leveraged to compute the right and left pseudoinverses efficiently. These are commonly computed for sparse systems using the singular-value-decomposition (SVD) \cite{Ben-Israel2003}, successive over-relaxation \cite{Björck1979}, and QR-factorization \cite{courrieu2008fast} methods. Note that $\bar{\mathbf{A}}^{\mathrm{T}}_{1}\bar{\mathbf{A}}_{1}$ is a non-singular\footnote{Note that $\bar{\mathbf{A}}_{1}$ has full column rank because it corresponds to non-zero components of a vertex of a polyhedron (see Proposition 2.1.4 of \citet{Bertsekas2009}).} sparse symmetric matrix, so we compute its inverse---computational bottleneck of pseudoinverses---using sparse LDL factorization \cite{ashcraft1998accurate}. Other matrix multiplications are performed using standard sparse operations. Second, because $\bar{\mathbf{A}}^{\mathrm{T}}_{1}$ has full row rank, an orthonormal basis of its null space can be reliably constructed using iterative SVD techniques for sparse matrices \cite{Berry1992}. Note that the computation of $\mathbf{Z}$ is not necessary in this section as it plays no role in the construction of $\Omega^{(k)}$. However, as will be discussed in Section~\ref{sec:multiple-d}, null-space computations are required when handling multiplicity using quadratic auxiliary objective.  

\subsection{Multiple Primal Solutions} \label{sec:multiple}
Next, we turn to the handling of primal multiplicity. We introduced two techniques in Section~\ref{sec:mpp-qp} to select a point from the optimal face of Eq.~(\ref{eqn:eq9}). In the first, continuous parametric solutions for $r$ decision variables are achieved by imposing a prescribed priority order through a hierarchical optimization problem with $r$ auxiliary linear objectives as follows
\begin{equation}%--------------------------------------------------------------------
		\left\{\begin{array}{l}
         	  z^{*}(\pmb{\theta})=\min\limits_{\mathbf{x}\in\mathbbm{R}^n} -\mathbf{c}^{\mathrm{T}}\mathbf{x}\\
						\mathbf{A}\mathbf{x}=\mathbf{w}+\mathbf{F}\pmb{\theta}\\
						\mathbf{x}\geq\mathbf{0}\\
						\left\{\begin{array}{l}
							z^{**}(\pmb{\theta})=\min\limits_{\mathbf{x}\in\mathbbm{R}^n} \mathbf{d}^{\mathrm{T}}_1\mathbf{x}\\
							\mathbf{A}\mathbf{x}=\mathbf{w}+\mathbf{F}\pmb{\theta}\\
							\mathbf{c}^{\mathrm{T}}\mathbf{x}=-z^{*}(\pmb{\theta})\\
							\mathbf{x}\geq\mathbf{0}\\
							\qquad\vdots
						\end{array} \right.
		\end{array} \right.,
	 \label{eqn:eq24}
\end{equation}%--------------------------------------------------------------------
where $\{\mathbf{d}_j\}_1^r$ are the auxiliary cost vectors. These must satisfy additional restrictions to ensure that Eq.~(\ref{eqn:eq24}) is well-posed. First, all $\mathbf{d}_j$ and rows of $\mathbf{A}$ must be linearly independent; otherwise, the entire optimal face of the primary LP remains optimal with respect to auxiliary objectives. Second, the order and sense of $\mathbf{d}_j$ must be appropriately chosen to avoid unboundedness in auxiliary LPs. Provided these conditions are satisfied, Eq.~(\ref{eqn:eq24}) can be solved with as many auxiliary objectives as necessary at $\pmb{\theta}=\pmb{\theta}_0$ to identify a vertex solution $\mathbf{x}^{*}(\pmb{\theta}_0)$ with $\mathcal{J}_1\equiv\{j\in\mathcal{J}|x^{*}_{j}(\pmb{\theta}_0)>0\}$ and $\mathcal{J}_0\equiv\{j\in\mathcal{J}|x^{*}_{j}(\pmb{\theta}_0)=0\}$. When $r<n$, the $r$th level solution may still not be unique; however, one may treat it as unique and proceed with the steps outlined in Sections~\ref{sec:unique-nd} and \ref{sec:unique-d} depending on the degeneracy of the vertex solution. This is a modeling decision one can make to ensure continuity for a subset of decision variables. Note that there always exists a vertex solution of the primary LP that is also a solution to all lower-level LPs (Theorem 2 of \citet{Harwood2015}). Simplex-based algorithms from the literature can be used to efficiently compute such vertex solutions (Algorithm 2 of \citet{Harwood2015} and \citet{Jones2007}). Another approach of note for handling hierarchical optimization is the equivalent-weight method \cite{Sherali1982, Sherali1983}. \citet{Sherali1983} proved the existence of a single-objective LP with equivalent weights, any optimal solution of which is also an optimal solution of a lexicographic LP; two algorithms were also provided for constructing the cost vector of the equivalent LP \cite{Sherali1982}. This cost vector furnishes the only auxiliary objective necessary to identify a vertex solution that satisfies the priority order.      

Two major challenges impede the application of the algorithms of \citet{Sherali1982} to metabolic networks. First, the equivalent weights grow exponentially with the number of hierarchical levels, causing integer overflow for large-scale problems. Second, finding a priority order that does not lead to unbounded optimal solutions, which often relies on exponential-complexity algorithms, is not straightforward. To avoid these issues, we propose an alternative algorithm in~\ref{sec:appendixc} to generate an equivalent cost vector for a single auxiliary LP with a unique optimal solution. The optimal solution of this LP remains unique in the entire parameter space, so the respective parametric solution is continuous for all decision variables. The proposed algorithm does not directly impose priority orders, enabling a fast and efficient computation of equivalent cost vectors.    

The lexicographic-LP approach outlined above is closely connected to the algorithm of \citet{Jones2007}. In the former, one constructs a continuous parametric solution with respect to all decision variables using $r=n$ auxiliary objectives according to an appropriate priority order of the decision variables. The latter imposes the priority order by weighting the decision variables in the primary objective through disparate perturbations of the cost-vector components. This amounts to a lexicographic-positivity check of two augmented matrices (see Definition 6 and Eq.~(9) of \citet{Jones2007}), constructed using the basis associated with each vertex of the optimal face. Accordingly, the optimal vertex that meets the priority order is identified once the necessary lexicographic-positivities are satisfied. 

The second technique was proposed by \citet{Spjotvold2007}, in which, using an auxiliary norm-minimizing objective, a unique (not necessarily a vertex) solution from the optimal face is identified. This strategy leads to a continuous parametric primal solution for all decision variables. Note that there is no special significance to this quadratic objective besides computational convenience. The continuity of the primal solution is guaranteed for any strictly convex objective function, including strictly convex quadratic objectives. This technique requires the index partition $\{\mathcal{J}_1,\mathcal{J}_0\}$, defined at the beginning of Section~\ref{sec:improvement}, which is not algorithmically straightforward to construct. As suggested in the literature \cite{Tondel2003, Spjotvold2005}, it can be determined using a non-vertex solution from interior-point algorithms. However, the resulting index partition can be highly unreliable since these algorithms always provide approximate solutions to within a tolerable duality gap. In the following, we propose an alternative technique to construct the index partition, requiring a basis of the optimal face. We choose $\mathbf{x}^{*}(\pmb{\theta}_0)$ to be the Chebyshev center of $X^{*}(\pmb{\theta}_0)$, a non-vertex solution maximally distanced from all the constraints, at which only the constraints in $\mathcal{A}^{*}(\pmb{\theta}_0)$ are active. Our approach relies on the concept of implicit inequalities, contrasting with that of \citet{Ben-Israel1968} based on $\{1\}$-generalized inverses.   

Let $X^{*}:=\{\mathbf{x}\in\mathbbm{R}^n|\mathbf{G}\mathbf{x}=\mathbf{v},\mathbf{x}\ge\mathbf{0}\}$ denote the optimal face of Eq.~(\ref{eqn:eq9}) at a given $\pmb{\theta}$, where $\mathbf{G}:=[\mathbf{A}^{\mathrm{T}}\: -\mathbf{c}]^{\mathrm{T}}$ and $\mathbf{v}:=[\mathbf{b}^{\mathrm{T}}\: z^{*}]^{\mathrm{T}}$ with $\mathbf{b}=\mathbf{w}+\mathbf{F}\pmb{\theta}$. The projection of $X^{*}$ onto the null space of $\mathbf{G}$ can be written 
\begin{equation}%--------------------------------------------------------------------
	X^{*}_{\mathcal{N}}=\{\mathbf{y}\in\mathcal{N}|\mathbf{N}\mathbf{y}+\mathbf{h}\ge\mathbf{0}\},
	 \label{eqn:eq25}
\end{equation}%--------------------------------------------------------------------
where $\mathbf{h}:=\mathbf{G}^{\mathrm{T}}(\mathbf{G}\mathbf{G}^{\mathrm{T}})^{-1}\mathbf{v}$, $\mathcal{N}:=\mathrm{ker}(\mathbf{G})$, $n_N:=\mathrm{dim}(\mathcal{N})=n-m-1$, and $\mathbf{N}$ is a matrix, the columns of which form a basis of $\mathcal{N}$. Note that the dimension of the optimal face $X^{*}$ can be controlled by varying the cost vector, while the dimension of $\mathcal{N}$ is fixed (Eq.~(\ref{eqn:eq9}) is assumed to remain well-conditioned). This implies that $\mathrm{dim}(X^{*}_{\mathcal{N}})\le n_N$, and some of the inequalities defining $X^{*}_{\mathcal{N}}$ can be implicit. Accordingly, Eq.~(\ref{eqn:eq25}) is rewritten 
\begin{equation}%--------------------------------------------------------------------
	X^{*}_{\mathcal{N}}=\{\mathbf{y}\in\mathcal{N}|\mathbf{N}_{\mathrm{ex}}\mathbf{y}+\mathbf{h}_{\mathrm{ex}}\ge\mathbf{0},\mathbf{N}_{\mathrm{im}}\mathbf{y}+\mathbf{h}_{\mathrm{im}}=\mathbf{0} \},
	 \label{eqn:eq26}
\end{equation}%--------------------------------------------------------------------
 where $\mathbf{N}_{\mathrm{ex}}$ and $\mathbf{N}_{\mathrm{im}}$ are sub-matrices of $\mathbf{N}$ corresponding to the explicit and implicit inequalities of Eq.~(\ref{eqn:eq25}); $\mathbf{h}_{\mathrm{ex}}$ and $\mathbf{h}_{\mathrm{im}}$ are similarly defined. Using $\mathbf{N}_{\mathrm{im}}$, one can decompose $\mathcal{N}$ into a null space and its orthogonal complement to obtain the smallest affine subspace of $\mathbbm{R}^n$ containing $X^{*}$. The projection of $X^{*}_{\mathcal{N}}$ onto the null space of $\mathbf{N}_{\mathrm{im}}$ is
\begin{equation}%--------------------------------------------------------------------
	X^{*}_{\mathcal{N},\mathcal{N}'}=\{\mathbf{u}\in\mathcal{N}'|\mathbf{N}_{\mathrm{ex}}\pmb{\aleph}\mathbf{u}+\mathbf{N}_{\mathrm{ex}}\mathbf{y}_{p}+\mathbf{h}_{\mathrm{ex}}\ge\mathbf{0}\},
	 \label{eqn:eq27}
\end{equation}%--------------------------------------------------------------------
where $\mathbf{y}_{p}:=-\mathbf{N}^{\dagger}_{\mathrm{im}}\mathbf{h}_{\mathrm{im}}$, $\mathcal{N}':=\mathrm{ker}(\mathbf{N}_{\mathrm{im}})$, $n_X:=\mathrm{dim}(\mathcal{N}')=n_N-\mathrm{rank}(\mathbf{N}_{\mathrm{im}})$, $\mathbf{N}^{\dagger}_{\mathrm{im}}$ is the right pseudoinverse of $\mathbf{N}_{\mathrm{im}}$, and $\pmb{\aleph}$ is a matrix, the columns of which form a basis of $\mathcal{N}'$. Note that $\mathcal{N}'$ is isomorphic to the smallest affine subspace containing $X^{*}$, so $\mathrm{dim}(X^{*})=n_X$. The Chebyshev center $\mathbf{u}_c$ of $X^{*}$ in $\mathcal{N}'$ can be obtained from Eq.~(\ref{eqn:eq27}) by solving an LP (see Section 8.5.1 of \citet{Boyd2004}). Accordingly, the Chebyshev center of $X^{*}$ with respect to the original decision variables is 
\begin{equation}%--------------------------------------------------------------------
	\mathbf{x}_c=\mathbf{N}(\pmb{\aleph}\mathbf{u}_c+\mathbf{y}_{p})+\mathbf{h}.
	 \label{eqn:eq28}
\end{equation}%--------------------------------------------------------------------
Implicit inequalities in Eq.~(\ref{eqn:eq26}) can be readily obtained using Theorem 3.2 of \citet{Telgen1982}, only requiring a solution of an auxiliary LP. Moreover, $\mathbf{N}_{\mathrm{im}}$ is generally rank-deficient. Therefore, pseudoinverse computations using the techniques discussed in Section~\ref{sec:unique-d} are not possible. Here, we directly compute the particular solution $\mathbf{y}_{p}$ using the least-square algorithm of \citet{Paige1982} for sparse systems. We compute an orthonormal basis of $\mathcal{N}'$ using approximate LU-based techniques \cite{Gotsman2008} as they perform more reliably than SVD methods for large-scale rank-deficient systems.  

In the next two sections, we demonstrate how to implement the strategy of \citet{Spjotvold2007} in our improved algorithm. The main idea is to partition the CR defined from the optimal active set of the primary LP Eq.~(\ref{eqn:eq10}) into several sub-CRs using the following auxiliary QP
\begin{equation}%--------------------------------------------------------------------
		\left\{\begin{array}{l}
         	  y^{*}(\pmb{\theta})=\min\limits_{\bar{\mathbf{x}}_1 \in\mathbbm{R}^{n_1}} \frac{1}{2}\bar{\mathbf{x}}^{\mathrm{T}}_1\bar{\mathbf{x}}_1\\
						\bar{\mathbf{A}}_1\bar{\mathbf{x}}_1=\mathbf{w}+\mathbf{F}\pmb{\theta}\\
						\bar{\mathbf{x}}_1\geq\mathbf{0}
		\end{array} \right.,
	 \label{eqn:eq29}
\end{equation}%--------------------------------------------------------------------
so that the resulting parametric primal solution in each sub-CR has the minimum distance from the origin among the optimal solution set of the primary LP. Equation~(\ref{eqn:eq29}) is to be solved at $\pmb{\theta}=\pmb{\theta}_0$ to ascertain the unique $\bar{\mathbf{x}}^{*}_1(\pmb{\theta}_0)$, which, in turn, furnishes a partition $\{\bar{\mathcal{J}}_{11},\bar{\mathcal{J}}_{10}\}$ of $\bar{\mathcal{J}}_1$ with $\bar{\mathcal{J}}_{11}:=\{j\in\bar{\mathcal{J}}_1|\bar{x}^{*}_{1,j}(\pmb{\theta}_0)>0\}$, $\bar{\mathcal{J}}_{10}:=\{j\in\bar{\mathcal{J}}_1|\bar{x}^{*}_{1,j}(\pmb{\theta}_0)=0\}$, $n_{11}:=|\bar{\mathcal{J}}_{11}|$, and $n_{10}:=|\bar{\mathcal{J}}_{10}|$. As with Eq.~(\ref{eqn:eqa1}), a secondary index map $\phi:\bar{\mathcal{J}}_1\longmapsto\bar{\bar{\mathcal{J}}}_1$ with $\Phi:=\phi^{-1}$ and its associated permutation matrix $\pmb{\mathcal{P}}$ are defined. The permutation matrix, viewed as a linear transformation $\pmb{\mathcal{P}}:\bar{\mathbf{x}}_1\mapsto\bar{\bar{\mathbf{x}}}_1$, separates all zero and non-zero variables of $\bar{\mathbf{x}}_1$ and their respective columns in the constraint matrix $\bar{\mathbf{A}}_1$. The forgoing map provides the index set $\bar{\bar{\mathcal{J}}}_1=\phi(\bar{\mathcal{J}}_1)$ and partition $\{\bar{\bar{\mathcal{J}}}_{11},\bar{\bar{\mathcal{J}}}_{10}\}$, decision-variable partition $(\bar{\bar{\mathbf{x}}}_{11},\bar{\bar{\mathbf{x}}}_{10})$, and dual-variable partition $(\bar{\bar{\pmb{\mu}}}_{1},\bar{\bar{\pmb{\mu}}}_{0})$ in the new coordinate system with $\bar{\bar{\mathbf{x}}}_{11}:=\bar{\bar{\mathbf{x}}}_{1,\bar{\bar{\mathcal{J}}}_{11}}=\bar{\mathbf{x}}_{1,\Phi(\bar{\bar{\mathcal{J}}}_{11})}$, $\bar{\bar{\mathbf{x}}}_{10}:=\bar{\bar{\mathbf{x}}}_{1,\bar{\bar{\mathcal{J}}}_{10}}=\bar{\mathbf{x}}_{1,\Phi(\bar{\bar{\mathcal{J}}}_{10})}$, $\bar{\bar{\pmb{\mu}}}_{1}:=\bar{\bar{\pmb{\mu}}}_{\bar{\bar{\mathcal{J}}}_{11}}$, and $\bar{\bar{\pmb{\mu}}}_{0}:=\bar{\bar{\pmb{\mu}}}_{\bar{\bar{\mathcal{J}}}_{10}}$. The transformed constraint matrix is $\bar{\bar{\mathbf{A}}}_1=[\bar{\bar{\mathbf{A}}}_{11}|\bar{\bar{\mathbf{A}}}_{10}]$, where $\bar{\bar{\mathbf{A}}}_{11}$ and $\bar{\bar{\mathbf{A}}}_{10}$ are sub-matrices containing all the columns of $\bar{\mathbf{A}}_1$ that correspond to non-zero and zero decision variables in $\bar{\mathbf{x}}_1$. Using these notations, Eq.~(\ref{eqn:eq29}) is rewritten with respect to the double-barred coordinate system
\begin{equation}%--------------------------------------------------------------------
		\left\{\begin{array}{l}
         	  y^{*}(\pmb{\theta})=\min\limits_{\bar{\bar{\mathbf{x}}}_1\in\mathbbm{R}^{n_1}} \frac{1}{2}\bar{\bar{\mathbf{x}}}^{\mathrm{T}}_1\bar{\bar{\mathbf{x}}}_1\\
						\bar{\bar{\mathbf{A}}}_1\bar{\bar{\mathbf{x}}}_1=\mathbf{w}+\mathbf{F}\pmb{\theta}\\
						\bar{\bar{\mathbf{x}}}_1\geq\mathbf{0}
		\end{array} \right.
	 \label{eqn:eq30}
\end{equation}%--------------------------------------------------------------------
with the respective KKT system
\begin{align}%--------------------------------------------------------------------
	\bar{\bar{\mathbf{x}}}_{11}-\bar{\bar{\mathbf{A}}}^{\mathrm{T}}_{11}\bar{\bar{\pmb{\lambda}}}&=\mathbf{0}_{n_{11}\times 1}, \label{eqn:eq31}\\
	\bar{\bar{\mathbf{A}}}^{\mathrm{T}}_{10}\bar{\bar{\pmb{\lambda}}}+\bar{\bar{\pmb{\mu}}}_{0}&=\mathbf{0}_{n_{10}\times 1}, \label{eqn:eq32}\\
	\bar{\bar{\mathbf{A}}}_{11}\bar{\bar{\mathbf{x}}}_{11}&=\mathbf{w}+\mathbf{F}\pmb{\theta},	\label{eqn:eq33}\\
	\bar{\bar{\mathbf{x}}}_{11}&\ge\mathbf{0}_{n_{11}\times 1},	\label{eqn:eq34}\\
	\bar{\bar{\pmb{\mu}}}_0&\ge\mathbf{0}_{n_{10}\times 1},	\label{eqn:eq35}
\end{align}%--------------------------------------------------------------------
with $\bar{\bar{\mathbf{x}}}_{10}$ and $\bar{\bar{\pmb{\mu}}}_1$ eliminated as they are identically zero in each sub-CR. Now, we have arrived at a suitable formulation to construct the solution triple $(\bar{\bar{\mathbf{x}}}^{*}_{11},\bar{\bar{\pmb{\lambda}}}^{*},\bar{\bar{\pmb{\mu}}}^{*}_0)$.

\subsubsection{Non-Degenerate Case} \label{sec:multiple-nd}
The minimum-norm solution $\bar{\bar{\mathbf{x}}}^{*}_{11}$ can be a vertex or non-vertex point, so this case applies to problems with $n_{11}\ge m$. Since the primal problem is non-degenerate, primal and dual solutions are all unique. When $n_{11}>m$ (non-vertex solution), Eq.~(\ref{eqn:eq33}) is an underdetermined system of equations for $\bar{\bar{\mathbf{x}}}^{*}_{11}$. Here, Eqs.~(\ref{eqn:eq31})-(\ref{eqn:eq33}) are coupled, so primal and dual solutions cannot be generally obtained from individual equations. However, noting that the right pseudoinverse provides the minimum-norm component of the solution in underdetermined systems, the primal solution can be directly obtained from Eq.~(\ref{eqn:eq33})
\begin{equation}%--------------------------------------------------------------------
		\bar{\bar{\mathbf{x}}}^{*}_{11}=\bar{\bar{\mathbf{A}}}^{\mathrm{T}}_{11}(\bar{\bar{\mathbf{A}}}_{11}\bar{\bar{\mathbf{A}}}^{\mathrm{T}}_{11})^{-1}(\mathbf{w}+\mathbf{F}\pmb{\theta}).
	 \label{eqn:eq36}
\end{equation}%--------------------------------------------------------------------
Note that this is possible only because of the particular choice of the quadratic objective in Eq.~(\ref{eqn:eq29}). For any other convex quadratic objective, it is necessary to decompose $\bar{\bar{\mathbf{x}}}^{*}_{11}$ into a minimum-norm and null-space component, similarly to how $\bar{\pmb{\lambda}}^{*}$ was treated in Section~\ref{sec:unique-d}. Moreover, the general solution of Eq.~(\ref{eqn:eq36}) applies to cases with $n_{11}=m$ (vertex solution), where $\bar{\bar{\mathbf{A}}}_{11}$ is a square matrix, because the Moore-Penrose inverse reduces to the standard non-singular inverse for balanced systems \cite{Ben-Israel2003}. The dual solutions can be constructed accordingly
\begin{align}%--------------------------------------------------------------------
	\bar{\bar{\pmb{\lambda}}}^{*}=(\bar{\bar{\mathbf{A}}}_{11}\bar{\bar{\mathbf{A}}}^{\mathrm{T}}_{11})^{-1}(\mathbf{w}+\mathbf{F}\pmb{\theta}),	\label{eqn:eq37}\\
	\bar{\bar{\pmb{\mu}}}^{*}_0=-\bar{\bar{\mathbf{A}}}_{10}(\bar{\bar{\mathbf{A}}}_{11}\bar{\bar{\mathbf{A}}}^{\mathrm{T}}_{11})^{-1}(\mathbf{w}+\mathbf{F}\pmb{\theta}).	
	\label{eqn:eq38}
\end{align}%--------------------------------------------------------------------
Unlike previous cases, there is no parameter-independent relationship between primal and dual feasibility for QPs, so they both define the boundaries of sub-CRs\footnote{In search Strategy I, sub-CRs need not be distinguished from other CRs. Due to the recursive nature of this strategy, other sub-CRs with the same optimal active set are visited in subsequent cycles of the algorithm.} 
\begin{multline}%--------------------------------------------------------------------
	\Omega^{(k)}=\{\pmb{\theta}\in R^{(k)}|\bar{\bar{\mathbf{A}}}^{\mathrm{T}}_{11}(\bar{\bar{\mathbf{A}}}_{11}\bar{\bar{\mathbf{A}}}^{\mathrm{T}}_{11})^{-1}(\mathbf{w}+\mathbf{F}\pmb{\theta})\ge\mathbf{0},\\
	\bar{\bar{\mathbf{A}}}_{10}(\bar{\bar{\mathbf{A}}}_{11}\bar{\bar{\mathbf{A}}}^{\mathrm{T}}_{11})^{-1}(\mathbf{w}+\mathbf{F}\pmb{\theta})\le\mathbf{0}\}. 
	 \label{eqn:eq39}
\end{multline}%--------------------------------------------------------------------

\subsubsection{Degenerate Case} \label{sec:multiple-d}
This case applies when $n_{11}<m$, where the optimal primal solution is a degenerate vertex. As a result, the dual variables have multiple optimal solutions, which we handle by the decomposition $\bar{\bar{\pmb{\lambda}}}^{*}=\pmb{\mathcal{Y}}\bar{\bar{\pmb{\lambda}}}^{*}_{y}+\pmb{\mathcal{Z}}\bar{\bar{\pmb{\lambda}}}^{*}_{z}$. The columns of $\pmb{\mathcal{Z}}$ and $\pmb{\mathcal{Y}}$ span the null space and its orthogonal complement of $\bar{\bar{\mathbf{A}}}^{\mathrm{T}}_{11}$. Following the same procedure as outlined in Section~\ref{sec:unique-d}, the primal and dual solutions are obtained 
\begin{align}%--------------------------------------------------------------------
\bar{\bar{\mathbf{x}}}^{*}_{11}&=(\bar{\bar{\mathbf{A}}}^{\mathrm{T}}_{11}\bar{\bar{\mathbf{A}}}_{11})^{-1}\bar{\bar{\mathbf{A}}}^{\mathrm{T}}_{11}(\mathbf{w}+\mathbf{F}\pmb{\theta}), \label{eqn:eq40}\\
\bar{\bar{\pmb{\lambda}}}^{*}&=\bar{\bar{\mathbf{A}}}_{11}(\bar{\bar{\mathbf{A}}}^{\mathrm{T}}_{11}\bar{\bar{\mathbf{A}}}_{11})^{-2}\bar{\bar{\mathbf{A}}}^{\mathrm{T}}_{11}(\mathbf{w}+\mathbf{F}\pmb{\theta})+\pmb{\mathcal{Z}}\bar{\bar{\pmb{\lambda}}}^{*}_{z}, \label{eqn:eq41}\\
	\bar{\bar{\pmb{\mu}}}^{*}_0&=-\bar{\bar{\mathbf{A}}}^{\mathrm{T}}_{10}\bar{\bar{\mathbf{A}}}_{11}(\bar{\bar{\mathbf{A}}}^{\mathrm{T}}_{11}\bar{\bar{\mathbf{A}}}_{11})^{-2}\bar{\bar{\mathbf{A}}}^{\mathrm{T}}_{11}(\mathbf{w}+\mathbf{F}\pmb{\theta})-\bar{\bar{\mathbf{A}}}^{\mathrm{T}}_{10}\pmb{\mathcal{Z}}\bar{\bar{\pmb{\lambda}}}^{*}_{z}. \label{eqn:eq42}
\end{align}%--------------------------------------------------------------------
As in the previous case, dual feasibility must be explicitly incorporated into the definition of the CR. The solutions in Eqs.~(\ref{eqn:eq40})-(\ref{eqn:eq42}) remain optimal provided there is a $\bar{\bar{\pmb{\lambda}}}^{*}_{z}\in\mathbbm{R}^{(m-n_{11})}$ such that Eqs.~(\ref{eqn:eq34}) and (\ref{eqn:eq35}) are satisfied. This can be formally stated as 
\begin{multline}%--------------------------------------------------------------------
	\Omega^{(k)}=\{\pmb{\theta}\in R^{(k)}|\exists\bar{\bar{\pmb{\lambda}}}^{*}_{z}\in\mathbbm{R}^{(m-n_{11})},(\bar{\bar{\mathbf{A}}}^{\mathrm{T}}_{11}\bar{\bar{\mathbf{A}}}_{11})^{-1}\bar{\bar{\mathbf{A}}}^{\mathrm{T}}_{11}(\mathbf{w}+\mathbf{F}\pmb{\theta})\ge\mathbf{0},\\
	\bar{\bar{\mathbf{A}}}^{\mathrm{T}}_{10}\bar{\bar{\mathbf{A}}}_{11}(\bar{\bar{\mathbf{A}}}^{\mathrm{T}}_{11}\bar{\bar{\mathbf{A}}}_{11})^{-2}\bar{\bar{\mathbf{A}}}^{\mathrm{T}}_{11}(\mathbf{w}+\mathbf{F}\pmb{\theta})+\bar{\bar{\mathbf{A}}}^{\mathrm{T}}_{10}\pmb{\mathcal{Z}}\bar{\bar{\pmb{\lambda}}}^{*}_{z}\le\mathbf{0}\}. 
	 \label{eqn:eq43}
\end{multline}%--------------------------------------------------------------------
Alternatively, $\Omega^{(k)}$ can be more conveniently represented as the projection of a $(q+m-n_{11})$-dimensional polyhedron $\hat{\Omega}^{(k)}$ onto a $q$-dimensional parameter space. $\hat{\Omega}^{(k)}$ is defined as  
\begin{multline}%--------------------------------------------------------------------
	\hat{\Omega}^{(k)}:=\{(\pmb{\theta},\bar{\bar{\pmb{\lambda}}}^{*}_{z})\in R^{(k)}\times\mathbbm{R}^{(m-n_{11})}|(\bar{\bar{\mathbf{A}}}^{\mathrm{T}}_{11}\bar{\bar{\mathbf{A}}}_{11})^{-1}\bar{\bar{\mathbf{A}}}^{\mathrm{T}}_{11}(\mathbf{w}+\mathbf{F}\pmb{\theta})\ge\mathbf{0},\\
	\bar{\bar{\mathbf{A}}}^{\mathrm{T}}_{10}\bar{\bar{\mathbf{A}}}_{11}(\bar{\bar{\mathbf{A}}}^{\mathrm{T}}_{11}\bar{\bar{\mathbf{A}}}_{11})^{-2}\bar{\bar{\mathbf{A}}}^{\mathrm{T}}_{11}(\mathbf{w}+\mathbf{F}\pmb{\theta})-\bar{\bar{\mathbf{A}}}^{\mathrm{T}}_{10}\pmb{\mathcal{Z}}\bar{\bar{\pmb{\lambda}}}^{*}_{z}\le\mathbf{0}\}, 
	 \label{eqn:eq44}
\end{multline}%--------------------------------------------------------------------
from which it follows that
\begin{equation}%--------------------------------------------------------------------
		\Omega^{(k)}=\mathrm{proj}_{R^{(k)}}\hat{\Omega}^{(k)}.
	 \label{eqn:eq45}
\end{equation}%--------------------------------------------------------------------
Equation~(\ref{eqn:eq45}) underpins the differences between auxiliary linear and quadratic objectives in handling primal multiplicity. Unlike LPs, for which dual feasibility is ensured by the strong duality theorem, irrespective of $\pmb{\theta}$, both primal and dual feasibility conditions must be satisfied at KKT points of QPs. This significantly complicates the construction of CRs for highly-degenerate primal problems ($m-n_{11}\gg1$). Since the optimal face of the dual for these problems is high-dimensional, several one-dimensional projection steps are required to construct the H-representation of $\Omega^{(k)}$ from $\hat{\Omega}^{(k)}$ algorithmically. Therefore, there is a trade-off between the convenience of achieving continuous parametric solutions and complexity of constructing CRs for quadratic and linear auxiliary objectives. 

Having derived a parametric primal solution with respect to transformed coordinate systems, we introduce additional notations to construct inverse transformations and determine the primal solution in the original coordinate system. Let $\mathbf{P}=[\mathbf{P}_1|\mathbf{P}_0]$, where $\mathbf{P}_1$ and $\mathbf{P}_0$ are sub-matrices containing the first $n_1$ and last $n_0$ columns of the permutation matrix used in Sections~\ref{sec:unique-nd} and ~\ref{sec:unique-d}. Similarly, let $\pmb{\mathcal{P}}=[\pmb{\mathcal{P}}_1|\pmb{\mathcal{P}}_0]$, where $\pmb{\mathcal{P}}_1$ and $\pmb{\mathcal{P}}_0$ are sub-matrices containing the first $n_{11}$ and last $n_{10}$ columns of the permutation matrix used in Sections~\ref{sec:multiple-nd} and ~\ref{sec:multiple-d}. Then, we have $\mathbf{x}^{*}=\mathbf{P}_1\bar{\mathbf{x}}^{*}_{1}$ and $\mathbf{x}^{*}=\mathbf{P}_1\pmb{\mathcal{P}}_1\bar{\bar{\mathbf{x}}}^{*}_{11}$ when the primal solution is unique and non-unique , respectively. Furthermore, many of the inequalities defining the CRs in Eqs.~(\ref{eqn:eq19}), (\ref{eqn:eq23}), (\ref{eqn:eq39}), and (\ref{eqn:eq43}) are redundant in most practical problems. To achieve a minimal H-representation, Theorem 2.2 of \citet{Telgen1982} can be used as it provides an algorithmically convenient criterion for identifying redundant inequalities from a system of linear constraints, only requiring a solution of an auxiliary LP. 
  
\section{Computational Complexity} \label{sec:complexity}
We highlighted the computational challenges of handling primal multiplicity using quadratic auxiliary objectives in Section~\ref{sec:multiple-d}. Difficulties arise when the optimal solution of this auxiliary QP is degenerate, where there are multiple optimal dual solutions. Here, a CR is constructed by projecting a $(q+d)$-dimensional polyhedron, defined by $h$ hyperplanes, onto a $q$-dimensional parameter space, where $d$ is the degeneracy degree of the primal solution. Although polyhedral-projection algorithms have significantly improved from the early Fourier-Motzkin method, they still suffer from undesirable scaling properties \cite{Tiwary2008}. This is due to many unnecessary constraints being generated in each one-dimensional projection, leading to a double-exponential complexity $O(h^{2^d})$ \cite{Monniaux2010}. Unnecessary constraints can be eliminated in each step using auxiliary LPs \cite{Telgen1982}, but the resulting single-exponential complexity is still not satisfactory \cite{Monniaux2010}.  

\begin{figure}
	\centering
	\includegraphics[width=\linewidth]{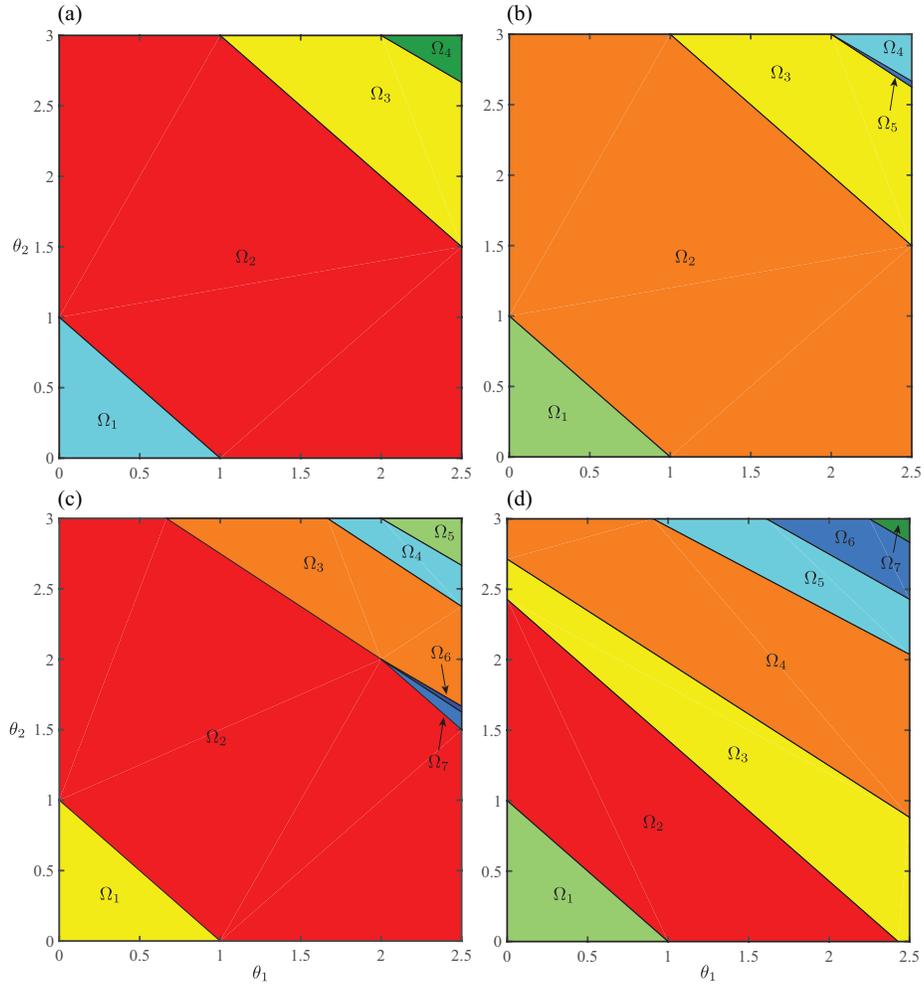}
  \caption{Polyhedral partition of the parameter space for Example 1 using (a) no auxiliary objective, (b) a linear auxiliary objective with $\mathbf{d}_1=[0 \: 0 \: -1]^{\mathrm{T}}$, (c) two linear auxiliary objectives with $\mathbf{d}_1=[0 \: -1 \: 0]^{\mathrm{T}}$ and $\mathbf{d}_2=[0 \: 0 \: 1]^{\mathrm{T}}$, and (d) an auxiliary norm-minimizing quadratic objective.}
	\label{fig:figure5}
\end{figure}

Similar challenges have to be overcome when handling degenerate optimal solutions. Eliminating redundant constraints from the set of active constraints is one approach. Given a system of linearly dependent active constraints $\mathbf{G}\mathbf{x}=\mathbf{b}$ with $\mathbf{x}\in\mathbbm{R}^n$, $\mathbf{G}\in\mathbbm{R}^{m\times n}$, $\mathbf{b}\in\mathbbm{R}^m$, and $m>n$, one needs to check the linear dependence of $(m-1)!$ constraint pairs to identify redundant equalities. Note that checking linear dependence of two vectors in $\mathbbm{R}^n$ requires at least $n$ arithmetic operations, so the overall complexity of this approach is estimated $O(m!)$. Alternatively, using QR factorization, a balanced system from the set of active constraints can be derived in polynomial time \cite{Borrelli2003}. Here, the Hessenberg factorization---one of the most efficient QR algorithms with a polynomial complexity $O(m^3)$---offers desirable scaling properties \cite{Arbenz2012}. Nevertheless, we handled degeneracies more conveniently using generalized inverses in Sections~\ref{sec:unique-d} and \ref{sec:multiple-d}. Here, the set of active constraints furnishes an overdetermined system of equations, so the solution is determined by the left pseudoinverse $\mathbf{x}=(\mathbf{G}^{\mathrm{T}}\mathbf{G})^{-1}\mathbf{G}^{\mathrm{T}}\mathbf{b}$ without needing to eliminate redundant equalities. The complexity of computing the left pseudoinverse relative to that of a balanced system with $n$ unknowns can readily be shown, by counting vector operations, to scale as $O(m^3)$. In addition, it can be implemented using sparse operations. Note that the $n\times n$ reduced matrices resulting from the first two approaches need to be inverted to provide the solution. Therefore, their respective complexities should also be interpreted as estimates relative to a balanced system of size $n$.   

\section{Numerical Examples} \label{sec:example}
Three examples are presented in this section, illustrating the implementation of our improved MPP algorithm. The first is taken from the work of \citet{Spjotvold2005}, where we compare the effectiveness of linear and quadratic auxiliary objectives in generating continuous parametric solutions. The second illustrates the application of MPP algorithms to a metabolic network model, where we discuss the physical implications of the solution characteristics. The third demonstrates the computational efficiency of our algorithm in handling large-scale metabolic networks.

\begin{figure} [ht]
	\centering
	\includegraphics[width=0.97\linewidth]{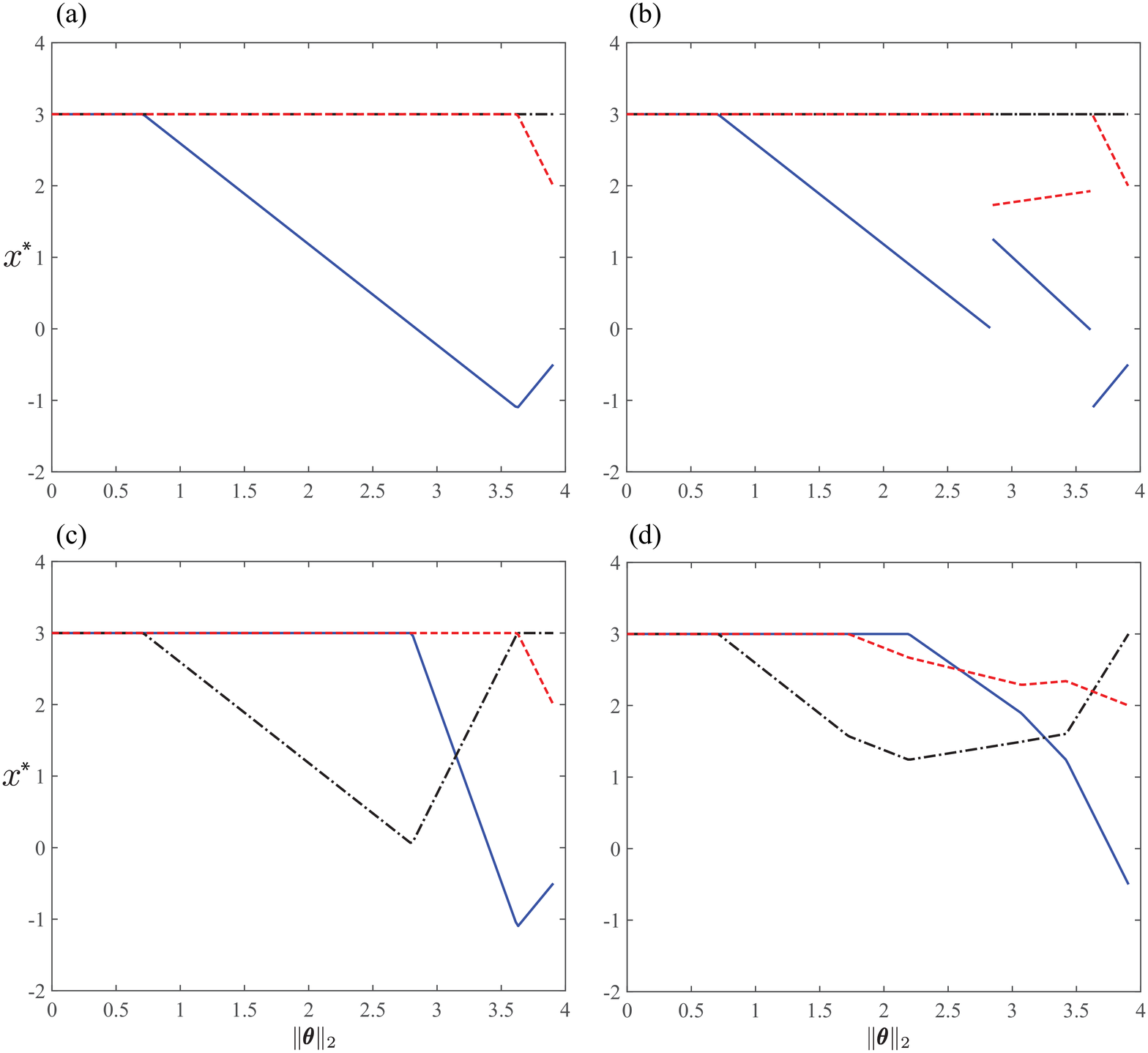}
  \caption{Continuity and discontinuity of the optimal primal solution $x^{*}_1$ (solid), $x^{*}_2$ (dashed), and $x^{*}_3$ (dashed-dotted) along the path $\mathcal{C}=\{\pmb{\theta}\in[0,2.5]\times[0,3]|2.5\theta_2=3\theta_1\}$ for Example 1 using the same auxiliary objectives as Fig.~\ref{fig:figure5}. The norm is measured from the lower left corner of the parameter region $R^{(1)}$.}
	\label{fig:figure6}
\end{figure}

\noindent
\textbf{Example 1.} Consider the parametric LP
\begin{equation}%--------------------------------------------------------------------
		\left\{\begin{array}{l}
         	  z^{*}(\pmb{\theta})=\max\limits_{\mathbf{x}\in\mathbbm{R}^{3}} \quad x_1+x_2+x_3\\
						x_1+x_2+x_3\le 10-\theta_1-\theta_2\\
						x_1-2x_2\le 4-\theta_1-2\theta_2\\
						-x_1-2x_3\le 3-\theta_1-2\theta_2\\
						|x_i|\le 3, \quad i\in\{1,2,3\}
		\end{array} \right.,
	 \label{eqn:eq46}
\end{equation}%--------------------------------------------------------------------
where a parametric solution in the region $R^{(1)}=[0,2.5]\times[0,3]$ is sought. This problem has multiple optimal primal solutions for all $\pmb{\theta}\in R^{(1)}$. Figure~\ref{fig:figure5} shows polyhedral partitions of $R^{(1)}$ using various linear and quadratic auxiliary objectives with their respective primal solutions along the path $2.5\theta_2=3\theta_1$ plotted in Fig.~\ref{fig:figure6}. Figure~\ref{fig:figure5}a (to be compared with Fig.~1a of \citet{Spjotvold2005}) shows a partition constructed without any auxiliary objective by arbitrarily selecting a vertex from the optimal face of Eq.~(\ref{eqn:eq46}). Although the optimal active sets in $\Omega_{1-4}$ are identical to those in the corresponding CRs in Fig.~1a of \citet{Spjotvold2005}, choosing a different optimal vertex led to a different partition of $R^{(1)}$. The resulting primal solution is continuous for our choice (Fig.~\ref{fig:figure6}a) and discontinuous for that of \citet{Spjotvold2005}. Using an auxiliary linear objective with $\mathbf{d}_1=[0 \quad 0 \quad -1]^{\mathrm{T}}$, resulting in the partition shown in Fig.~\ref{fig:figure5}b, ensures a continuous solution only for $x^{*}_3$ (Fig.~\ref{fig:figure6}b). Continuous solutions for all decision variables can be achieved using two linear auxiliary objectives with $\mathbf{d}_1=[0 \quad -1 \quad 0]^{\mathrm{T}}$ and $\mathbf{d}_2=[0 \quad 0 \quad 1]^{\mathrm{T}}$ (Fig.~\ref{fig:figure6}c) at the expense of generating more CRs (Fig.~\ref{fig:figure5}c).

Figure~\ref{fig:figure5}d (to be compared with Fig.~1b of \citet{Spjotvold2005}) shows a partition constructed using an auxiliary norm-minimizing objective. As expected, the resulting primal solution is continuous for all decision variables (Fig.~\ref{fig:figure6}d). Note that this partition is different than that of \citet{Spjotvold2005} because of the standard-form formulation we adopted: firstly, unlike the LP representation of Example 1, a sign change in each decision variable triggers an active-set change in standard form. Consequently, applying MPP algorithms to standard-form problems generally results in more CRs than to problems with unrestricted decision variables. Secondly, all artificial and slack variables contribute to the norm function of the auxiliary QP in our formulation (see Eq.~(\ref{eqn:eq29})), whereas \citet{Spjotvold2005} only considered the norm of the original decision variables in Example 1. 

\begin{figure} 
	\centering
	\includegraphics[width=\linewidth]{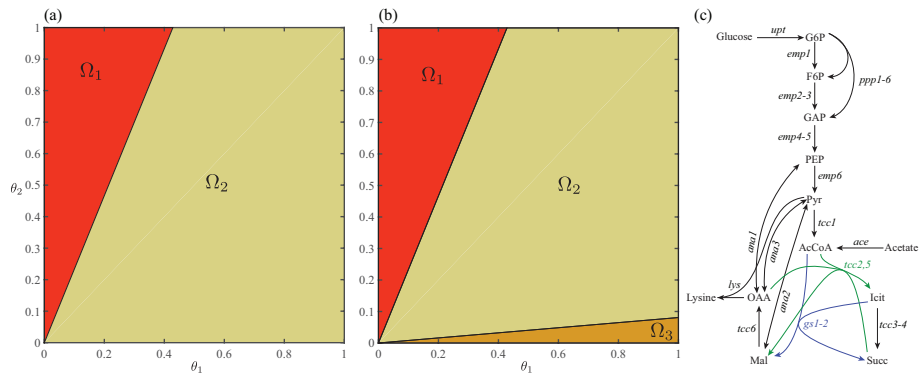}
  \caption{Polyhedral partition of the parameter space for a reduced metabolic model of lysine biosynthesis \cite{Schuster2007} in Example 2 using (a) no auxiliary objective and (b) an auxiliary norm-minimizing objective. Substrate uptake rates are bounded by Michaelis-Menten kinetic rates \cite{Dowd1965}. Here, $\theta_1=-v^{\mathrm{lb}}_g/v_{\mathrm{max},g}$ and $\theta_2=-v^{\mathrm{lb}}_a/v_{\mathrm{max},a}$ measure the maximum uptake rates of glucose and acetate scaled with their respective maximum initial velocities $v_{\mathrm{max},g}=2.1$ \cite{Lindner2011} and $v_{\mathrm{max},a}=1.8$ \cite{Ebbighausen1991} (mmol/gDW/h). (c) Schematic of lysine network adapted from \citet{Schuster2007}. Reaction labels denote major enzyme-catalyzed pathways.}
	\label{fig:figure7}
\end{figure}

Metabolic-network models based on flux balance analysis (FBA) conform to the MPP formulation presented in Section~\ref{sec:improvement}. These usually comprise a complex network of reactions, through which reactants are converted to products. Their primary goal is to identify reaction rates and the optimal pathway a micro-organism follows to produce its main product, namely, biomass \cite{Schilling1998}. FBA simplifies complicated dynamic metabolic models by assuming a steady-state intracellular metabolism, postulating that optimal pathways maximize the growth (biomass production) rate, thereby, avoiding direct computations of reaction rates from kinetic constants \cite{Orth2010}. Reaction rates are the unknowns of FBA, constrained by the network stoichiometry and represented by sign-restricted (irreversible reactions) and sign-unrestricted (reversible reactions) variables; reactions with zero rates indicate inactive\footnote{Not to be confused with the inactive set. In the standard-form formulation of FBA, the inactive set corresponds to non-zero decision variables, which are associated with active reactions in the network. A similar statement can be made about the active set.} parts of the network. Therefore, FBA can be formulated as an LP, in which the growth rate is maximized subject to linear stoichiometric constraints. The following example demonstrates how MPP can be applied to study metabolic networks using FBA. 

\noindent
\textbf{Example 2.} The metabolic-network model of lysine biosynthesis from glucose and acetate in \emph{C. glutamicum} (Fig.~\ref{fig:figure7}c) \cite{Schuster2007} can be formulated as the LP
\begin{equation}%--------------------------------------------------------------------
		\left\{\begin{array}{l}
         	  z^{*}(\pmb{\theta})=\max\limits_{\mathbf{v}\in\mathbbm{R}^{n'}} v_l\\
						\mathbf{S}\mathbf{v}=\mathbf{0}\\
						v_j\ge v_j^{\mathrm{lb}}, \quad j\in\{g,a\}\\
						v_j\ge 0, \quad j\in\mathcal{J}_{\mathrm{irr}}
		\end{array} \right.,
	 \label{eqn:eq47}
\end{equation}%--------------------------------------------------------------------
where $\mathbf{v}$ is the vector of reaction rates, $\pmb{\theta}=-[v^{\mathrm{lb}}_g/v_{\mathrm{max},g} \: v^{\mathrm{lb}}_a/v_{\mathrm{max},a}]^{\mathrm{T}}$, and $\mathcal{J}_{\mathrm{irr}}$ is the index set of irreversible reactions with $a$, $g$, and $l$ the indexes of acetate uptake, glucose uptake, and lysine production reactions. The network stoichiometry matrix $\mathbf{S}\in\mathbbm{R}^{m'\times n'}$ is determined from the reaction equations listed in Table 1 of \citet{Schuster2007}. This model comprises $m'=33$ equality constraints (metabolites) and $n'=35$ decision variables (reactions) in its non-standard form Eq.~(\ref{eqn:eq47}). As with many other species \cite{Dowd1965}, the maximum substrate uptake rates of \emph{C. glutamicum} follow Michaelis-Menten kinetics \cite{Ebbighausen1991, Lindner2011}
\begin{equation}%--------------------------------------------------------------------
	 v_j^{\mathrm{lb}}=\frac{-v_{\mathrm{max},j}C_j}{K_{\mathrm{m},j}+C_j}, \quad j\in\{g,a\},\\
	 \label{eqn:eq48}
\end{equation}%--------------------------------------------------------------------
where maximum initial velocities $v_{\mathrm{max},j}$ and saturation constants $K_{\mathrm{m},j}$ are the kinetic parameters and $C_j$ are the substrate concentrations in the extracellular environment. Maximum uptake rates furnish the lower bounds in Eq.~(\ref{eqn:eq47}) because reaction rates corresponding to the consumption of metabolites from the extracellular environment are negative by convention in FBA \cite{Schilling1998}. Hence, the parametric solution of Eq.~(\ref{eqn:eq47}) for $\pmb{\theta}\in[0,1]^2$ describes the metabolism of \emph{C. glutamicum} with respect to every possible concentrations of glucose and acetate in its growth environment. 

Polyhedral partitions for the LP in Eq.~(\ref{eqn:eq47}) with and without auxiliary objectives are plotted in Fig.~\ref{fig:figure7}. Here, we emphasize two important features of MPP algorithms in relation to FBA and their physical implications. First, the algorithm we presented accounts for decision-variable signs. Recall, CRs represent regions of the parameter space where an active set remains optimal. Since active sets in standard form correspond to zero decision variables, the inactive set of each CR in FBA can be regarded as representing a distinct metabolic pathway. Accordingly, crossing CR boundaries in the parameter space corresponds to either a major alteration of the optimal pathway or a minor direction change in active reversible reactions. Therefore, there are intimate connections between CRs and metabolic modes of an organism, reflecting its metabolic adaptations and regulatory mechanisms in response to external stimuli. These metabolic modes are comparable to the elementary modes arising from elementary mode analysis \cite{Stelling2002}. The standard-form formulation particularly serves this interpretation by distinguishing between reversible and irreversible reactions \cite{Schuster2000} through sign-sensitive active sets. For example, consider the partition constructed using no auxiliary objective (Fig.~\ref{fig:figure7}a) with two CRs, the active sets of which differ only in NAD and NADH reaction rates. The production rate of NAD (NADH) is positive (negative) in $\Omega_1$ and negative (positive) in $\Omega_2$. Hence, $\Omega_1$ represents glucose-deprived states associated with low-yield lysine production using oxidative phosphorylation. In contrast, $\Omega_2$ represents acetate-deprived states corresponding to high-yield lysine production without using oxidative phosphorylation \cite{Schuster2007}.    

The physical interpretation of the optimal active set is also of interest when multiple solutions exit. The optimal active set of a CR represents pathways that are inactive for every possible solution of the optimal-solution set. For the lysine network in Fig.~\ref{fig:figure7}c, the optimal active set includes the reactions \emph{tcc1}, \emph{tcc3}, \emph{tcc4}, and \emph{tcc5} from the tricarboxylic acid cycle for all the CRs in Figs.~\ref{fig:figure7}a and b. Therefore, these reactions can never contribute toward maximum-lysine production. In the context of elementary mode analysis, this implies that the elementary modes comprising these reactions play no role in the optimal pathway, irrespective of glucose and acetate concentrations.      

Second, the continuity of parametric solutions---a natural regularity condition---imposed by applying auxiliary objective functions in our algorithm is motivated by its application to metabolic networks, given that organisms are not expected to abruptly alter their metabolism in response to environmental stimuli. Figure~\ref{fig:figure7}b shows a polyhedral partition of Eq.~(\ref{eqn:eq47}) associated with a continuous parametric solution for all physically feasible parameters. It is constructed using an auxiliary norm-minimizing objective, which splits the second CR of the primary LP ($\Omega_2$ in Fig.~\ref{fig:figure7}a) into two sub-CRs ($\Omega_2$ and $\Omega_3$ in Fig.~\ref{fig:figure7}b), leaving the first CR unchanged. The metabolic modes associated with $\Omega_1$ in Fig.~\ref{fig:figure7}a and Fig.~\ref{fig:figure7}b are the same, so are those associated with $\Omega_2$ in Fig.~\ref{fig:figure7}a and its respective sub-CRs $\Omega_2$ and $\Omega_3$ in Fig.~\ref{fig:figure7}b. The optimal pathways corresponding to these sub-CRs differ only in the reaction \emph{emp6} (see Fig.~\ref{fig:figure7}c), which is inactive in the former and active in the latter. A careful examination of empirical data at low acetate concentrations is necessary to determine whether this minor change in the optimal pathway is of physical relevance. 

We conclude this example by highlighting the physical significance of the facets arising from partitioning the parameter space. Consider the facets $\partial \Omega_{1,2}$ (Type-I hyperplane) and $\partial \Omega_{2,3}$ (Type-II hyperplane) in Fig.~\ref{fig:figure7}b. The first signifies a transition from a glucose-deprivation to acetate-deprivation induced metabolism in the parameter space, whereas the second results from a norm-minimizing auxiliary objective applied to ensure a continuous parametric solution. This regularity condition can be imposed by other convex objectives, but the resulting facet and optimal pathways may not be the same. Note that auxiliary objectives may be applied for reasons other than imposing continuity. As has been recognized in the literature \cite{Stelling2002}, FBA cannot capture all cellular regulation using maximal growth alone. Additional constraints may be imposed on the optimal solution space of metabolic models through other auxiliary objectives to achieve more realistic predictions of cell functions over a wide range of environmental conditions.

\begin{figure} [t]
	\centering
	\includegraphics[width=\linewidth]{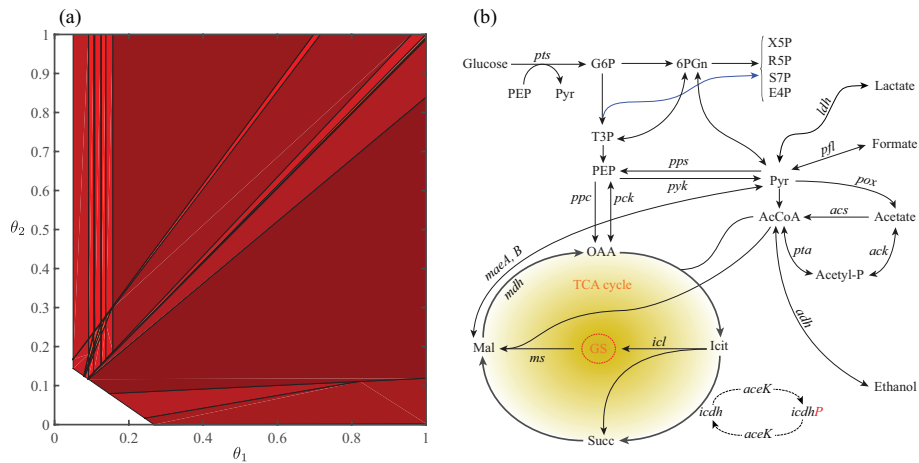}
  \caption{(a) Polyhedral partition of the parameter space for a genome-scale model of \emph{E. coli} (iJR904) \cite{Reed2003} in Example 3 using no auxiliary objective. Here, $\theta_1=-v^{\mathrm{lb}}_g/v_{\mathrm{max},g}$ and $\theta_2=-v^{\mathrm{lb}}_o/v_{\mathrm{max},o}$ are the maximum uptake rates of glucose and oxygen scaled with their respective maximum initial velocities $v_{\mathrm{max},g}=10.5$ and $v_{\mathrm{max},o}=15$ (mmol/gDW/h) \cite{Hanly2011}. Xylose metabolism is excluded from the network by setting the respective upper and lower bounds to zero. The unshaded area represents an infeasible region of the respective LP problem, where a steady-state metabolism cannot be sustained. (b) Compressed model of \emph{E. coli} metabolism adapted from \citet{Castano-Cerezo2009}. Reaction labels denote the enzymes catalyzing major pathways.}
	\label{fig:figure8}
\end{figure}

\noindent
\textbf{Example 3.} A metabolic-network model of \emph{E. coli} (iJR904) \cite{Reed2003} is studied, emphasizing the computational challenges facing MPP algorithms when handling genome-scale metabolic networks (see Fig.~\ref{fig:figure8}b). Physical interpretations of the results will be discussed elsewhere. Here, the metabolism is described with respect to the glucose and oxygen concentrations through the LP model  
\begin{equation}%--------------------------------------------------------------------
		\left\{\begin{array}{l}
         	  z^{*}(\pmb{\theta})=\max\limits_{\mathbf{v}\in\mathbbm{R}^{n'}} v_x\\
						\mathbf{S}\mathbf{v}=\mathbf{0}\\
						v_j\ge v_j^{\mathrm{lb}}, \quad j\in\{g,o\}\\
						v_j\ge v_j^{\mathrm{lb,fix}}, \quad j\in\mathcal{J}_{\mathrm{lb,fix}}\\
						v_j\le v_j^{\mathrm{ub,fix}}, \quad j\in\mathcal{J}_{\mathrm{ub,fix}}\\
						v_j\ge 0, \quad j\in\mathcal{J}_{\mathrm{irr}}	
		\end{array} \right.,
	 \label{eqn:eq49}
\end{equation}%--------------------------------------------------------------------
where $\pmb{\theta}=-[v^{\mathrm{lb}}_g/v_{\mathrm{max},g} \: v^{\mathrm{lb}}_o/v_{\mathrm{max},o}]^{\mathrm{T}}$, $\mathcal{J}_{\mathrm{lb,fix}}$ is the index set of reactions with a fixed lower bound, $\mathcal{J}_{\mathrm{ub,fix}}$ is the index set of reactions with a fixed upper bound, and $\mathcal{J}_{\mathrm{irr}}$ is the index set of irreversible reactions with $o$, $g$, and $x$ the indices of oxygen uptake, glucose uptake, and biomass production reactions, respectively. This model comprises $m'=761$ metabolites and $n'=1075$ reactions. The xylose metabolism is excluded from the network to limit the number of parameters. Accordingly, the xylose uptake rate is fixed at zero by imposing the respective upper and lower bounds. Fixed upper and lower bounds in Eq.~(\ref{eqn:eq49}) include the foregoing uptake-rate bounds and a lower bound on the ATP exchange reaction, accounting for ATP maintenance requirements \cite{Reed2003}. 

\begin{figure} [t]
	\centering
	\includegraphics[width=\linewidth]{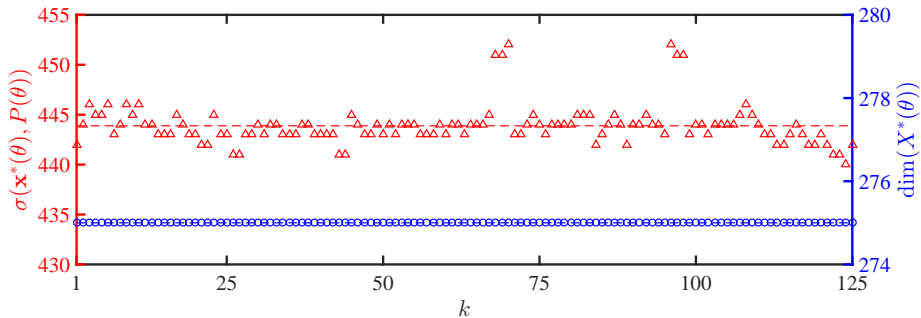}
  \caption{Optimal-solution characteristics of Example 3 for all CRs depicted in Fig.~\ref{fig:figure8}. The degeneracy degree (triangles) at a vertex solution $\mathbf{x}^{*}$ of the feasible set $P$ and the dimension (circles) of the optimal face $X^{*}$ with their respective averages (dashed lines) are plotted versus the CR counter $k$ for a $\theta\in\Omega^{(k)}$.}
	\label{fig:figure9}
\end{figure}

A polyhedral partition of the parameter region $R^{(1)}=[0,1]^2$ for the LP in Eq.~(\ref{eqn:eq49}) using the maximal-growth objective is depicted in Fig.~\ref{fig:figure8}a. The dimension of the optimal face and degeneracy degree at a primal-solution vertex for the respective CRs are shown in Fig.~\ref{fig:figure9}. Evidently, this LP has highly degenerate ($\sigma\approx 450$) and multiple optimal solutions with respect to its primary objective for all $\pmb{\theta}\in R^{(1)}$. Therefore, resolving multiplicity using a quadratic auxiliary objective, as discussed in Section~\ref{sec:multiple-d}, is impractical. Moreover, all optimal faces are high-dimensional, suggesting that the construction of continuous parametric solutions using any multiplicity-resolving algorithm is computationally challenging. Note also that there are several CRs sharing a facet with more than one CR. Thus, Strategy II cannot be applied using the algorithms of \citet{Tondel2003a} and \citet{Jones2007}. Among available techniques, the algorithm of \citet{Jones2006} and the lexicographic-LP approach introduced in Section~\ref{sec:multiple} are tractable for constructing continuous parametric solutions of large-scale metabolic networks. 

\begin{figure} [t]
	\centering
	\includegraphics[width=\linewidth]{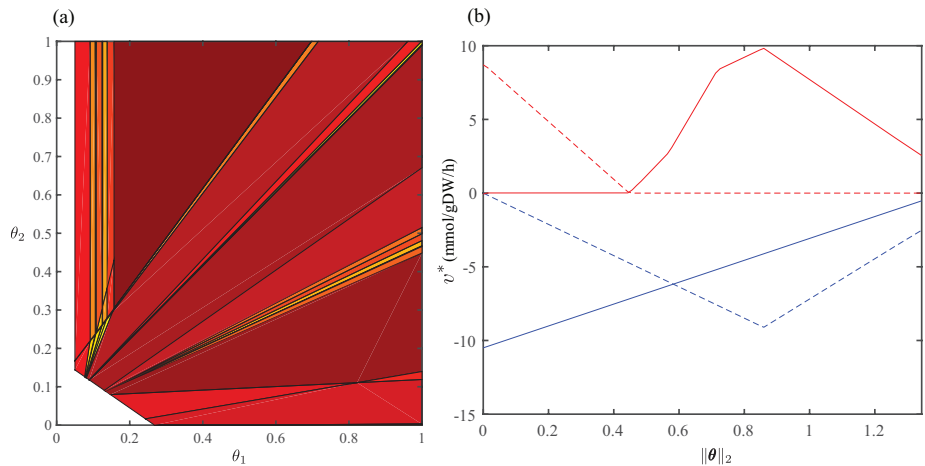}
  \caption{(a) Polyhedral partition of Example 3 constructed using the same parameters as for Fig.~\ref{fig:figure8} and an auxiliary objective, the equivalent cost vector of which furnishes continuous parametric solution with respect to all decision variables (see \ref{sec:appendixc}). (b) Continuity of optimal reaction rates for important metabolites, including glucose (solid blue), oxygen (dashed blue), carbon dioxide (solid red), and ethanol (dashed red) along the path $\mathcal{C}=\{\pmb{\theta}\in[0,1]^2|\theta_2=1-\theta_1\}$, where the norm is measured from the lower right corner of the parameter region $R^{(1)}$.}
	\label{fig:figure10}
\end{figure}

Figure~\ref{fig:figure10}a shows a polyhedral partition of the LP in Eq.~(\ref{eqn:eq49}) with 71 CRs constructed using maximal growth as the primary objective and an auxiliary objective with an equivalent cost vector generated by the algorithm outlined in~\ref{sec:appendixc}. The multiplicity criterion of \citet{Appa2002} was checked in all CRs to ensure that the optimal solution of the auxiliary LP is unique. Consequently, the corresponding parametric solution is continuous for all decision variables in $R^{(1)}$, as illustrated in Fig.~\ref{fig:figure10}b for the optimal reaction rates for key metabolites. Note that, in general, infinitely many cost vectors can be equivalent to the same lexicographic objective that imposes a given priority order. Moreover, many feasible priority orders (\ie, those having bounded optimal solution at all levels) may exist for a given LP. Therefore, one can generally find several polyhedral partitions of the parameter space for which the parametric solution is continuous. 

\section{Computational Performance} \label{sec:performance}
In this section, we compare the computational performance of our algorithm with that of Multi-Parametric Toolbox (MPT) \cite{Herceg2013}---an open source software package for parametric optimization. All studied cases were run in MATLAB R2016a on a 64-bit Windows machine with a 3.07 GHz Intel Xeon CPU using the LP solver CPLEX. We examine four cases of varying size, including the three examples studied in Section~\ref{sec:example} and a reduced version of iJR904, namely, \emph{E. coli} core \cite{Orth2010a}. Computation times for our algorithm are reported for cases where primal multiplicity is resolved using a single auxiliary LP with an equivalent cost vector constructed using the technique discussed in \ref{sec:appendixc}. Parametric solutions furnished by our algorithm and MPT are continuous with respect to all decision variables. 

\begin{table}[t]
\caption{Computational performance of the present algorithm versus Multi-Parametric Toolbox (MPT) \cite{Herceg2013}. The computation time of the equivalent cost vector $t_c$, partitioning the parameter space $t_p$, and their total $t$ are reported separately for the present algorithm. The number of equality constraints $m$ and decision variables $n$ in standard form indicate problem size.}
\centering
\begin{tabular}{l @{\hskip.3cm} c @{\hskip.1cm} c @{\hskip.5cm} c @{\hskip-.7cm} c @{\hskip-.8cm} c @{\hskip.0cm} c @{\hskip.5cm} c @{\hskip.1cm} c}
\hline
\shortstack[l]{} & \shortstack[c]{}  & \shortstack[c]{} & \shortstack[c]{} & \shortstack[c]{\hskip.4cmPresent algorithm} & \shortstack[c]{}  & \shortstack[c]{} & \shortstack[c]{} & \shortstack[c]{\hskip-1cm MPT} \\
\hline
\shortstack[l]{} & \shortstack[c]{$m$} & \shortstack[c]{$n$} & \shortstack[c]{$t_c$ (s)}  & \shortstack[c]{$t_p$ (s)} & \shortstack[c]{$t$ (s)} & \shortstack[c]{\#CRs} & \shortstack[c]{$t$ (s)}  & \shortstack[c]{\#CRs} \\
\hline
J.S.$^{\dagger}$ & 9 & 15 & 0.013 & 2.136 & 2.149 & 5 & 1.731 & 7  \\
\emph{C. glutamicum} & 33 & 55 & 0.015 & 1.926 & 1.941 & 2 & 2.139 & 2  \\
\emph{E. coli} core & 72 & 144 & 0.042 & 11.024 & 11.066 & 19 & 6.203 & 19 \\
\emph{E. coli} iJR904 & 761 & 1337 & 63.490 & 60.447 & 123.937 & 71 & 570.848 & 51  \\
\hline
\end{tabular}
\flushleft
\small
$^{\dagger}$ Parametric LP in Eq.~(\ref{eqn:eq46}) adapted from the work of \citet{Spjotvold2005} studied in Example 1.
\label{tbl:table2}
\end{table}

Table~\ref{tbl:table2} summarizes the computational performance of our algorithm versus MPT. One can observe that, for small-scale problems (\ie, $O(100)$ equality constraints and decision variables in standard form), the computational speedup and overhead of sparse linear algebra are of the same order, so that the overall performance is case dependent. However, for large-scale problems (\ie, $>O(1000)$ equality constraints and decision variables), the efficiency of sparse computations outweighs the overheads, enhancing the overall performance. This indicates the computational advantages of the proposed algorithm for genome-scale metabolic networks.  

\section{Conclusions} \label{sec:conclusion}
A multi-parametric programming algorithm based on active-set methods was presented, improving on the algorithms of \citet{Bemporad2002} and \citet{Borrelli2003} in handling primal degeneracy and multiplicity for large-scale sparse systems. This was motivated by application to flux balance analysis where metabolic networks are modeled as parametric linear programs. Flux balance analysis is a convenient tool, providing reaction rates and optimal metabolic pathways as functions of extracellular conditions (\eg, substrate concentrations) without relying on exact intracellular kinetic data. However, linear programs arising from metabolic models of most organisms have multiple and highly degenerate optimal solutions. For degenerate cases, constructing a balanced system of equations from the optimality conditions (\ie, KKT system) is challenging due to computational complexities associated with eliminating redundant constraints. When there are multiple optimal solutions, selecting the ``best'' among all possible solutions is not straightforward because, in ranges of environmental conditions that are of practical interest, not every choice can properly represent all regulatory mechanisms and cell functions.

To handle primal degeneracy, we implemented the Moore-Penrose generalized inverse in the formulation of the optimality conditions. This furnished an efficient technique to treat underdetermined and overdetermined systems arising from primal degeneracy by avoiding search mechanisms required for eliminating redundant constraints. Furthermore, algebraic structures arising from generalized inverses can be leveraged by several robust sparse-decomposition techniques, enabling fast and efficient computations of parametric solutions. To handle primal multiplicity, we examined two approaches. In the first, an auxiliary norm-minimizing objective was introduced into the primary linear program to identify a unique primal solution. This provides a continuous parametric solution for all decision variables, but construction of critical regions requires several one-dimensional polyhedral projections for highly degenerate problems---a computationally expensive step. In the second, primal multiplicity was resolved via a lexicographic linear program. This technique avoids the complexities of polyhedral projections; however, the continuity of parametric solutions is not guaranteed for all decision variables. Continuous parametric solutions can be achieved for more decision variables at higher computational cost by adding more auxiliary linear objectives. Depending on the application at hand, either technique can be more computationally desirable then the other.

A genome-scale metabolic model of \emph{E. coli} was studied, illustrating the computational challenges facing multi-parametric programming algorithms. The linear program arising from this model has a high-dimensional and highly degenerate optimal face for all extracellular concentrations of glucose and oxygen, rendering the foregoing quadratic-programming approach impractical. Handling multiplicity using an auxiliary lexicographic linear program, which can only access vertex solutions, is tractable. However, noting that vertices constitute a small subset of all admissible solutions to a metabolic model, vertex solutions may not satisfy all necessary biological constraints to accurately describe various cell functions. Tractable algorithms that can handle highly degenerate and multiple (vertex or non-vertex) optimal solutions for large-scale systems are yet to be developed. Nevertheless, compared to existing multi-parametric programming tools, the proposed algorithm proved more efficient in handling large-scale metabolic networks.    

Application of parametric programming to flux balance analysis of metabolic networks partitions their parameter space into critical regions, each can be associated with an optimal metabolic mode. Maximal growth is a widely accepted criterion for identifying optimal pathways. However, the maximal-growth objective alone rarely results in a unique partition for most organisms of interest. Here, a crucial question arises as to whether there is a partition that best describes all metabolic modes and cell functions. Imposing regularity conditions, such as continuity, through auxiliary objectives is one approach that we examined in this work. However, many auxiliary objectives can be used to ensure continuity, each of which may result in different partitions and optimal pathways. Whether there are particular objective functions that can universally characterize all metabolic networks or if additional biologically-motivated restrictions (\eg, network robustness \cite{Stelling2002}) and regularity conditions need to be imposed (through auxiliary objectives or other constraint types) remain open questions for future studies. 

\appendix
\section{Construction of Permutation Matrices} \label{sec:appendixa}
A permutation matrix $\mathbf{P}$ for rearranging the columns of a matrix $\mathbf{A}$ can be associated with a bijective index map (see Chapter 0.3 of \citet{Ben-Israel2003}) 
\begin{equation}%--------------------------------------------------------------------
	 \left\{\begin{array}{l}
         	  \xi:\mathcal{J}\longmapsto\bar{\mathcal{J}},\quad\Xi:=\xi^{-1}\\
						\bar{j}=\xi(j)\\
						j=\Xi(\bar{j})
		\end{array} \right..
	 \label{eqn:eqa1}
\end{equation}%--------------------------------------------------------------------
Such a permutation was introduced in Section~\ref{sec:improvement} for each CR to transform the original coordinate system so that all non-zero (zero) decision variables are collected at the top (bottom) of the primal solution vector in the transformed coordinate system $\bar{\mathbf{x}}^{*}$. Here, $\{\Xi(\bar{j})|\bar{j}=1,2,\cdots, n_1\}$ ($\{\Xi(\bar{j})|\bar{j}=n_1+1, n_1+2, \cdots, n\}$) can be thought of as an ordered set containing the indexes of non-zero (zero) decision variables in the original coordinate, where $\Xi(\bar{j})$ is readily ascertained from $\mathbf{x}^{*}(\pmb{\theta}_0)$. Accordingly, $\mathbf{P}$ is defined
\begin{equation}%--------------------------------------------------------------------
	 P_{lk}:=\left\{\begin{array}{l}
         	  0\qquad l\neq\Xi(k)\\
						1\qquad l=\Xi(k)
		\end{array} \right.,\quad\forall l\in\mathcal{J},\quad \forall k\in\bar{\mathcal{J}}.
	 \label{eqn:eqa2}
\end{equation}%--------------------------------------------------------------------
Transformation of the constraint matrix $\mathbf{A}$ under $\mathbf{P}$ furnishes a matrix $\bar{\mathbf{A}}$, in which all the columns corresponding to non-zero (zero) decision variables are collected on the left (right). A convenient representation of this transformation can be written 
\begin{equation}%--------------------------------------------------------------------
	 \left\{\begin{array}{l}
         	  \bar{\mathbf{A}}=\mathbf{A}\mathbf{P}\\
						\bar{\mathbf{x}}=\mathbf{P}^{\mathrm{T}}\mathbf{x}
		\end{array} \right.,
	 \label{eqn:eqa3}
\end{equation}%--------------------------------------------------------------------
where the inverse transformations can readily be obtained from the orthogonality of the permutation matrix $\mathbf{P}\mathbf{P}^{\mathrm{T}}=\mathbf{I}_n$. Any vector $\mathbf{v}\in\mathbbm{R}^n$ that is associated with the space of decision variables transforms similarly to $\mathbf{x}$.

\section{Minimum-Norm Dual Solution and Strong Duality} \label{sec:appendixb}
We shall demonstrate that there always exists a $\bar{\pmb{\lambda}}^{*}_{z}\in\mathbbm{R}^{(m-n_1)}$, for which dual feasibility Eq.~(\ref{eqn:eq15}) holds. Specifically, we show that $\bar{\pmb{\mu}}^{*}_0(\bar{\pmb{\lambda}}^{*}_{z}=\mathbf{0})$ is always a dual solution of the KKT system in Eqs.~(\ref{eqn:eq11})-(\ref{eqn:eq15}). Following standard procedures from linear programming \cite{Bertsimas1997}, we consider a feasible direction $\bar{\mathbf{u}}$ corresponding to an optimal solution $\bar{\mathbf{x}}^{*}$ of Eq.~(\ref{eqn:eq10}) such that $\bar{\mathbf{x}}^{*}+\bar{\mathbf{u}}\ge\mathbf{0}$ and $\bar{\mathbf{A}}(\bar{\mathbf{x}}^{*}+\bar{\mathbf{u}})=\mathbf{w}+\mathbf{F}\pmb{\theta}$. Let $(\bar{\mathbf{u}}_1,\bar{\mathbf{u}}_0)$ be a partition of $\bar{\mathbf{u}}$ corresponding to the permutation matrix defined in Section~\ref{sec:unique-d}. Note that $\bar{\mathbf{u}}_0\ge\mathbf{0}$ because zero variables can only increase along feasible directions. Since $\bar{\mathbf{A}}\bar{\mathbf{x}}^{*}=\mathbf{w}+\mathbf{F}\pmb{\theta}$, we have $\bar{\mathbf{A}}\bar{\mathbf{u}}=\mathbf{0}$, from which it follows 
\begin{equation}%--------------------------------------------------------------------
	\bar{\mathbf{A}}_1\bar{\mathbf{u}}_1=-\bar{\mathbf{A}}_0\bar{\mathbf{u}}_0.
	 \label{eqn:eqb1}
\end{equation}%--------------------------------------------------------------------
Equation~(\ref{eqn:eqb1}) provides an overdetermined system of equations for $\bar{\mathbf{u}}_1$ (see Section~\ref{sec:unique-d}), the solution of which is determined using the left pseudoinverse of $\bar{\mathbf{A}}_1$
\begin{equation}%--------------------------------------------------------------------
	\bar{\mathbf{u}}_1=-(\bar{\mathbf{A}}_1^{\mathrm{T}}\bar{\mathbf{A}}_1)^{-1}\bar{\mathbf{A}}_1^{\mathrm{T}}\bar{\mathbf{A}}_0\bar{\mathbf{u}}_0.
	 \label{eqn:eqb2}
\end{equation}%--------------------------------------------------------------------
Moreovere, since $\bar{\mathbf{x}}^{*}$ is an optimal solution, the objective value can only increase along any feasible direction. Therefore, $z(\bar{\mathbf{x}}^{*}+\bar{\mathbf{u}})-z(\bar{\mathbf{x}}^{*})\ge\mathbf{0}$, and, consequently, 
\begin{equation}%--------------------------------------------------------------------
	\bar{\mathbf{c}}_1^{\mathrm{T}}\bar{\mathbf{u}}_1+\bar{\mathbf{c}}_0^{\mathrm{T}}\bar{\mathbf{u}}_0\le\mathbf{0},
	 \label{eqn:eqb3}
\end{equation}%--------------------------------------------------------------------
from which it follows that 
\begin{equation}%--------------------------------------------------------------------
	\bar{\mathbf{c}}_0^{\mathrm{T}}-\bar{\mathbf{c}}_1^{\mathrm{T}}(\bar{\mathbf{A}}_1^{\mathrm{T}}\bar{\mathbf{A}}_1)^{-1}\bar{\mathbf{A}}_1^{\mathrm{T}}\bar{\mathbf{A}}_0\le\mathbf{0}.
	 \label{eqn:eqb4}
\end{equation}%--------------------------------------------------------------------
On the other hand, letting $\bar{\pmb{\lambda}}^{*}_{z}=\mathbf{0}$, the dual solution in Eq.~(\ref{eqn:eq22}) can be written 
\begin{equation}%--------------------------------------------------------------------
	(\bar{\pmb{\mu}}^{*}_0)^{\mathrm{T}}=-\bar{\mathbf{c}}_0^{\mathrm{T}}+\bar{\mathbf{c}}_1^{\mathrm{T}}(\bar{\mathbf{A}}_1^{\mathrm{T}}\bar{\mathbf{A}}_1)^{-1}\bar{\mathbf{A}}_1^{\mathrm{T}}\bar{\mathbf{A}}_0.
	 \label{eqn:eqb5}
\end{equation}%--------------------------------------------------------------------
Therefore, dual feasibility $\bar{\pmb{\mu}}^{*}_0\ge\mathbf{0}$ is implied by Eqs.~(\ref{eqn:eqb4}) and (\ref{eqn:eqb5}). This shows that the minimum-norm component of $\bar{\pmb{\lambda}}^{*}$ furnishes a feasible $\bar{\pmb{\mu}}^{*}_0$, which can also be interpreted as a consequence of the strong duality theorem for LPs. 

\section{Equivalent Cost Vectors} \label{sec:appendixc}
Consider the first-level auxiliary LP of Eq.~(\ref{eqn:eq24})
\begin{equation}%--------------------------------------------------------------------
	z^{**}=\min\limits_{\mathbf{x}\in\mathbbm{R}^n} \mathbf{d}^{\mathrm{T}}\mathbf{x}
							\quad\mbox{s.t.}\quad\mathbf{G}\mathbf{x}=\mathbf{h},
							\quad\mathbf{x}\geq\mathbf{0}
	 \label{eqn:eqc1}
\end{equation}%--------------------------------------------------------------------
with the respective dual 
\begin{equation}%--------------------------------------------------------------------
							y^{**}=\max\limits_{\pmb{\lambda}\in\mathbbm{R}^{m+1}} \mathbf{h}^{\mathrm{T}}\pmb{\lambda}
							\quad\mbox{s.t.}\quad\mathbf{G}^\mathrm{T}\pmb{\lambda}\le\mathbf{d},
	 \label{eqn:eqc2}
\end{equation}%--------------------------------------------------------------------
where $\mathbf{G}:=[\mathbf{A}^{\mathrm{T}}\: -\mathbf{c}]^{\mathrm{T}}$ and $\mathbf{h}:=[\mathbf{b}^{\mathrm{T}}\: z^{*}]^{\mathrm{T}}$ with $\mathbf{b}=\mathbf{w}+\mathbf{F}\pmb{\theta}$. We are concerned with ranges of $\pmb{\theta}$, for which the primary LP has multiple optimal solutions, so that Eq.~(\ref{eqn:eqc1}) is feasible. We seek to find a $\mathbf{d}\in\mathbbm{R}^n$ such that the optimal solution of Eq.~(\ref{eqn:eqc1}) is bounded and unique. Equation~(\ref{eqn:eqc1}) with such a cost vector is equivalent to a lexicographic LP that furnishes a unique solution by imposing an $n$-level priority order. Accordingly, we refer to $\mathbf{d}$ as an \emph{equivalent cost vector}. The technique presented here ensures the uniqueness and continuity of all decision variables of the primary LP as a whole rather than explicitly imposing them on individual variables.

The forgoing boundedness and uniqueness requirements can be imposed through the dual problem Eq.~(\ref{eqn:eqc2}) to establish a desirable range for $\mathbf{d}$. Since Eq.(\ref{eqn:eqc1}) is feasible by assumption, the boundedness of its optimal solution is guaranteed by the feasibility of Eq.(\ref{eqn:eqc2}). The feasible set of Eq.~(\ref{eqn:eqc2}), which can generally be an unbounded polyhedron, is bounded by an arbitrarily large positive number $U$ as
\begin{equation}%--------------------------------------------------------------------
	\Lambda(\mathbf{d}):=\{\pmb{\lambda}\in\mathbbm{R}^{m+1}:\mathbf{G}^\mathrm{T}\pmb{\lambda}\le\mathbf{d},|\pmb{\lambda}|\le\mathbf{1}U\}.
	 \label{eqn:eqc3}
\end{equation}%--------------------------------------------------------------------
The two inequalities in Eq.(\ref{eqn:eqc3}) can be combined as $\hat{\mathbf{G}}^\mathrm{T}\pmb{\lambda}\le\hat{\mathbf{d}}$, where $\hat{\mathbf{G}}:=[\mathbf{G}\:\:\mathbf{L}]$, $\hat{\mathbf{d}}:=[\mathbf{d}^\mathrm{T}\:\:\mathbf{1}_{1\times2(m+1)}U]^\mathrm{T}$, and $\mathbf{L}:=[\mathbf{I}_{m+1}\:\:-\mathbf{I}_{m+1}]$. In general, some of the $2(m+1)$ added inequalities $|\pmb{\lambda}|\le\mathbf{1}U$ are redundant and need not be included in Eq.(\ref{eqn:eqc3}). Assuming that $s$ out of $2(m+1)$ added inequalities are non-redundant, Eq.(\ref{eqn:eqc3}) can be represented more compactly as $\Lambda(\mathbf{d})=\{\pmb{\lambda}\in\mathbbm{R}^{m+1}:\mathbf{G}'^{\mathrm{T}}\pmb{\lambda}\le\mathbf{d}'\}$, where $\mathbf{G}':=[\mathbf{G}\:\:\mathbf{L}']$. Here, $\mathbf{L}'\in\mathbbm{R}^{(m+1)\times s}$ and $\mathbf{d}'\in\mathbbm{R}^{n+s}$ are, respectively, a sub-matrix of $\mathbf{L}$ and a sub-vector of $\hat{\mathbf{d}}$, in which redundant inequalities are eliminated. Let 
\begin{equation}%--------------------------------------------------------------------
	D:=\{\mathbf{d}\in\mathbbm{R}^{n}:\Lambda(\mathbf{d})\neq\emptyset\}
	 \label{eqn:eqc4}
\end{equation}%--------------------------------------------------------------------
be the set of desirable directions. We seek an explicit representation of $D$ as linear inequalities only with respect to $\mathbf{d}$. One can choose an appropriate measure of the center $\pmb{\lambda}_c(\mathbf{d})$ of $\Lambda$, so that $\tilde{D}:=\{\mathbf{d}\in\mathbbm{R}^{n}:\mathbf{G}'^{\mathrm{T}}\pmb{\lambda}_c(\mathbf{d})\le\mathbf{d}'\}$ is a reasonable approximation of $D$. The geometric center is an appropriate choice since it always lies inside convex sets and $\tilde{D}=D$. However, given an H-representations of a polyhedral set, no explicit expression for geometric center is available. Here, we choose $\pmb{\lambda}_c(\mathbf{d})$ to be the point of minimal Euclidean norm from all the facets of $\Lambda(\mathbf{d})$ because it can be computed conveniently using the left pseudoinverse $\pmb{\lambda}_c(\mathbf{d})=(\mathbf{G}'\mathbf{G}'^{\mathrm{T}})^{-1}\mathbf{G}'\mathbf{d}'$. This center, however, does not always lie inside the respective set, even for convex sets. Nevertheless, it provides an expedient approximation $\tilde{D}\subseteq D$ of the set of desirable directions. Introducing the partition $[\mathbf{Q}|\mathbf{W}]=(\mathbf{G}'\mathbf{G}'^{\mathrm{T}})^{-1}\mathbf{G}'$, where $\mathbf{Q}\in\mathbbm{R}^{(m+1)\times n}$ and $\mathbf{W}\in\mathbbm{R}^{(m+1)\times s}$, the inequalities defining $\tilde{D}$ are obtained 
\begin{align}%--------------------------------------------------------------------
	(\mathbf{G}^\mathrm{T}\mathbf{Q}-\mathbf{I}_n)\mathbf{d}&\le-U\mathbf{G}^\mathrm{T}\mathbf{W}\mathbf{1}_{s\times1}, \label{eqn:eqc5}\\
	\mathbf{L}'^{\mathrm{T}}\mathbf{Q}\mathbf{d}&\le U(\mathbf{I}_{s}-\mathbf{L}'^{\mathrm{T}}\mathbf{W})\mathbf{1}_{s\times1}.  \label{eqn:eqc6}
\end{align}%--------------------------------------------------------------------
In the limit $U\to\infty$, Eq.(\ref{eqn:eqc6}) reduces to $|d_i|<\infty$, while the right-hand side of Eq.(\ref{eqn:eqc5}) approaches $\mathbf{0}$, so that $(\mathbf{G}^\mathrm{T}\mathbf{Q}-\mathbf{I}_n)\mathbf{d}\le\mathbf{0}$ defines an approximate cone $\tilde{D}$ of desirable directions for $\Lambda(\mathbf{d})$.

To ensure primal uniqueness, a $\mathbf{d}\in\tilde{D}$ must be chosen such that the optimal solution of Eq.(\ref{eqn:eqc2}) is non-degenerate for $\pmb{\theta}$ in a parameter region of interest. In other words, one must avoid directions that are perpendicular to any $(n-k)$-dimensional faces at the intersection of any set of $k$ hyperplanes defining the primal feasible set. Hence, if a $\mathbf{d}\in\tilde{D}$ is randomly chosen, the corresponding optimal solution will be unique with probability one. To construct such a random direction, let $\tilde{D}'$ be the set of directions in $\tilde{D}$ that are bounded by an arbitrarily large positive number $U$
\begin{equation}%--------------------------------------------------------------------
	\tilde{D}':=\{\mathbf{d}\in\mathbbm{R}^{n}:(\mathbf{G}^\mathrm{T}\mathbf{Q}-\mathbf{I}_n)\mathbf{d}\le\mathbf{0},|\mathbf{d}|\le\mathbf{1}U\}
	 \label{eqn:eqc7}
\end{equation}%--------------------------------------------------------------------
with $r_c$ and $\mathbf{d}_c$ the respective Chebyshev radius and Chebyshev center. Since $\mathbf{d}_c$ is not a random vector, randomness is introduced into the choice of equivalent cost vector according to
\begin{equation}%--------------------------------------------------------------------
	\mathbf{d}=\mathbf{d}_c+A\mathbf{t},
	 \label{eqn:eqc8}
\end{equation}%--------------------------------------------------------------------
where $\mathbf{t}\in[-1,1]^n$ is a random vector with $A$ the respective amplitude. Noting that $\tilde{D}$ is generally a proper subset of $D$, a reasonable lower bound for the amplitude is $r_c/U$.

\section*{Acknowledgments}
This work was funded by the Cooperative Agreement between the Masdar Institute of Science and Technology (Masdar Institute), Abu Dhabi, UAE and the Massachusetts Institute of Technology (MIT), Cambridge, MA, USA - Reference 02/MI/MI/CP/11/07633/GEN/G/00 for work under the Second Five Year Agreement.

\section*{References}
\bibliography{mybibfile}

\end{document}